\documentclass[12pt]{amsart}
\usepackage{amsmath,amsthm,amsfonts,amssymb,eucal}
\renewcommand {\a}{ \alpha }
\renewcommand{\b}{\beta}
\newcommand{\e}{\epsilon}

\newcommand{\g}{\gamma}
\newcommand{\G}{\Gamma}
\newcommand{\U}{\Upsilon}

\newcommand{\vark}{\varkappa}
\renewcommand{\d}{\delta}
\newcommand{\s}{\sigma}
\renewcommand{\l}{\lambda}
\renewcommand{\L}{\Lambda}
\newcommand{\z}{\zeta}
\renewcommand{\t}{\theta}

\newcommand{\p}{\partial}
\newcommand{\om}{\omega}
\newcommand{\Om}{\Omega}

\newcommand{\oq}{\ {\raise 7pt\hbox{${\scriptstyle\circ}$}}
\kern -7pt{%\lower 2pt
\hbox{$Q$}}}

\newcommand{\R}{ \mathbb R}

\newcommand{\Rtwo}{ \mathbb R^2}

\newcommand{\Ttwo}{\mathbb T^2}
\newcommand{\Ztwo}{\mathbb Z^2}

\newcommand {\GD}{\mathfrak D}

\newcommand {\GH}{\mathfrak H}

\newcommand {\ba}{\mathbf a}
\newcommand {\boldf}{\mathbf f}
\newcommand {\boldg}{\mathbf g}

\newcommand {\bb}{\mathbf b}

\newcommand {\BA}{\mathbf A}
\newcommand {\BB}{\mathbf B}
\newcommand {\BD}{\mathbf D}
\newcommand {\BG}{\mathbf G}

\newcommand {\BF}{\mathbf F}
\newcommand {\BM}{\mathbf M}
\newcommand {\BJ}{\mathbf J}
\newcommand {\BI}{\mathbf I}

\newcommand {\BR}{\mathbf R}
\newcommand {\BT}{\mathbf T}

\newcommand {\bx}{\mathbf x}

\newcommand {\be}{\mathbf e}
\newcommand {\bh}{\mathbf h}

\newcommand {\bz}{\mathbf z}
\newcommand {\by}{\mathbf y}

\newcommand {\bn}{\mathbf n}

\newcommand {\brho}{\boldsymbol\rho}

\newcommand {\bxi}{\boldsymbol\xi}
\newcommand {\bg}{\boldsymbol\gamma}
\newcommand {\bsig}{\boldsymbol\sigma}
\newcommand {\BOLG}{\boldsymbol\Gamma}
\newcommand {\BXI}{\boldsymbol\Xi}
\newcommand {\BSIG}{\boldsymbol\Sigma}

\newcommand{\lu}{\langle}
\newcommand{\ru}{\rangle}

\newcommand{\CO}{\mathcal O}
\newcommand{\CP}{\mathcal P}

\newcommand{\CC}{\mathcal C}
\newcommand{\CS}{\mathcal S}
\newcommand{\CE}{\mathcal E}
\newcommand{\CD}{\mathcal D}

\newcommand{\boldW}[2]{ \textup{\textbf{\textsf{W}}}^{#1, #2}}
\newcommand{\plainW}[2]{\textup{{\textsf{W}}}^{#1, #2}}
\newcommand{\plainC}[1]{\textup{{\textsf{C}}}^{#1}}
\newcommand{\boldC}[1]{\textup{\textbf{\textsf{C}}}^{#1}}
\newcommand{\boldH}[1]{\textup{\textbf{\textsf{H}}}^{#1}}
\newcommand{\plainH}[1]{\textup{{\textsf{H}}}^{#1}}
\newcommand{\boldL}[1]{\textup{\textbf{\textsf{L}}}^{#1}}
\newcommand{\plainL}[1]{\textup{{\textsf{L}}}^{#1}}

 \DeclareMathOperator {\im
}{{Im}} \DeclareMathOperator {\re} {{Re}}

\DeclareMathOperator{\mes}{{mes}}

\DeclareMathOperator{\loc}{\textup{\tiny loc}}

 \DeclareMathOperator{\dc}{d}
\DeclareMathOperator{\Mu}{M}

\newtheorem{thm}{Theorem}[section]
\newtheorem{cor}[thm]{Corollary}
\newtheorem{lem}[thm]{Lemma}
\newtheorem{prop}[thm]{Proposition}
\newtheorem{cond}[thm]{Condition}

\theoremstyle{definition}
\newtheorem{defn}[thm]{Definition}%[section]

\theoremstyle{remark}
\newtheorem{rem}[thm]{Remark}

\numberwithin{equation}{section}

\newcommand{\bee}{\begin{equation}}
\newcommand{\ene}{\end{equation}}
\newcommand{\bes}{\begin{split}}
\newcommand{\ens}{\end{split}}

\newcommand{\bet}{\begin{tm}}
\newcommand{\ent}{\end{tm}}
\newcommand{\bel}{\begin{lm}}
\newcommand{\enl}{\end{lm}}
\newcommand{\bec}{\begin{cor}}
\newcommand{\enc}{\end{cor}}
\newcommand{\bep}{\begin{pr}}
\newcommand{\enp}{\end{pr}}
\newcommand{\ber}{\begin{rem}}
\newcommand{\enr}{\end{rem}}

\newcommand{\Z}{\mathbb Z}

\newcommand{\CF}{\mathcal F}

\makeatletter
\def\square{\RIfM@\bgroup\else$\bgroup\aftergroup$\fi
  \vcenter{\hrule\hbox{\vrule\@height.6em\kern.6em\vrule}\hrule}\egroup}
\makeatother

\begin{document}
\hoffset -4pc

\title[Quasi-conformal mappings]
{Quasi-conformal mappings
and periodic spectral problems
in dimension two}
\author[E. Shargorodsky and A.V. Sobolev
%(\the\day.\the\month.\the\year)
]
{Eugene Shargorodsky and Alexander V. Sobolev}
\address{Centre for Mathematical Analysis
and Its Applications\\
University of Sussex\\
Falmer,  Brighton\\
BN1 9QH, UK}
\email{E.Shargorodsky@sussex.ac.uk,
A.V.Sobolev@sussex.ac.uk}

\date{26 July 2001}
\begin{abstract}
We study spectral properties of second order elliptic operators
with periodic coefficients in dimension two.
These operators act in periodic
simply-connected waveguides,
with either Dirichlet, or Neumann, or the third
boundary condition. The main result is the absolute continuity
of the  spectra of such operators.
The corner stone of the proof
is an isothermal change of variables,
reducing the metric to a flat one and the
waveguide to a straight strip.
The main technical tool is the quasi-conformal variant of the
Riemann mapping theorem.
\end{abstract}

\maketitle

\section{Introduction}

According to a common belief, second order elliptic differential
operators with periodic coefficients should not have degenerate
bands in their spectra, or, in other words, their spectra should
be purely absolutely continuous (see \cite{BS}, \cite{KuchLev},
\cite{S}). The first rigorous proof of this fact was given by L.
Thomas in \cite{Th} for the Schr\"odinger operator $\Delta+V$ with
a periodic real-valued potential $V$. Further developments in this
area were driven by the attempts to consider operators with ever
increasing ``strength'' of the periodic perturbation, i.e. to pass
from zero order (as in\cite{Th}) to first and second order
perturbations, and ultimately to tackle the absolute continuity of
the elliptic operator
\begin{equation}\label{ham:eq}
H = \sum_{j, l =1}^d(D_j - a_j)g_{jl}(D_l - a_l) + V,
\ \
D_j = -i\p_j,
\end{equation}
with a periodic variable metric
$\{g_{jl}\} = \BG$ and a magnetic
potential $\ba = \{a_l\}$,
for arbitrary $d\ge 2$.

The case of first order perturbations,
i.e. that of constant $\BG$'s
and variable $\ba$'s,  was
handled  in \cite{HH} (small $\ba$'s),
\cite{BS1}, \cite{BS2} ($d=2$)
and \cite{Sob} (arbitrary $d\ge 2$).
If the metric is conformal, i.e. the matrix
$\BG$ is given by a scalar multiple of
the identity matrix, then the problem
can be easily reduced to the case of a
constant metric.
This situation is discussed in \cite{BS}.
The most difficult case, that  of a general variable
$\BG$, remained unaccessible until the paper
\cite{M}, where it was resolved for infinitely differentiable
$\BG$'s, $\ba$'s and $d = 2$.
An important breakthrough was made
in the recent work \cite{F}
where the absolute continuity was proved
in all dimensions $d\ge 2$, but with an additional
requirement of the reflectional symmetry of the operator.
Without the symmetry assumption the
question is still open. At present it is only known that
without smoothness assumptions on the coefficients
the absolute continuity may break down. An appropriate
example with $V = 0, \ba = 0$ and a ``non-smooth'' $\BG$
was constructed in \cite{Fil}.

Most of the progress was achieved in
the two-dimensional case, which we shall
discuss in more details.
The paper \cite{M}, as well as all the earlier
papers on absolute continuity,
relied on the approach suggested by L. Thomas in \cite{Th}.
Later it was observed in \cite{KuchLev} that
the general periodic metric can be reduced
to a conformal (or scalar) periodic one
by a suitable isothermal change of variables.
This allows one
to reproduce the results of \cite{M}
even under much weaker
assumptions, by reducing the
problem to the one considered in \cite{BS1},
\cite{BS2}, \cite{BS}. This approach was
exploited in \cite{BSS} for the operator of the form
\eqref{ham:eq} with a delta-like
periodic potential supported
by a periodic system of curves. Even more general perturbations
are studied in \cite{Sh1}, \cite{Sh2}.

The main aim of the present paper
is to prove the absolute continuity for
the operator of the form \eqref{ham:eq}
with non-constant periodic
coefficients,
defined in a periodic
domain $\Om\subset\Rtwo$
(which is usually referred to as a waveguide)
with the Dirichlet
 or ``natural'' boundary conditions.
In these two cases we use the
notation $H_D$ or $H_N$ for the
operator at hand.
By the natural boundary conditions we mean either
Neumann or the third boundary condition.
As in \cite{BSS}, we also include in the operator
a delta-like periodic perturbation supported by
smooth curves, which allows us to study
the cases of the Neumann condition and
the third boundary condition simultaneously.
In a somewhat restricted generality
this problem
with Dirichlet and Neumann conditions
was considered in \cite{SobWal}.
Compared to \cite{SobWal} in the present paper
the smoothness
conditions on the coefficients are substantially relaxed,
and the case of the third boundary-value problem is also
treated. The progress has become
possible due to a different approach to the problem:
instead of Morame's techniques we now use
the isothermal coordinate change.
A similar approach is used in the work
\cite{ShSus} focused on this problem in a
slightly different setting.

The corner stone of
our method is an isothermal change of
coordinates, which reduces the metric to a
conformal one.
The new coordinates are given by functions satisfying
the Beltrami equation (see Sect \ref{conf6:sect} below)
with a dilatation coefficient
$q$
determined
by the matrix $\BG$ (see \eqref{q:eq}),
i.e. they define a
$q$-quasi-conformal mapping.
We use a
$q$-quasi-conformal change of variables which
maps the waveguide $\Om$
homeomorphically onto a straight strip.
The existence of such a transformation
follows from a ``quasi-conformal'' version
of Riemann's mapping theorem (see e.g. \cite{Kru}, Ch. 1).
Here a key point is that
the uniqueness part of Riemann's mapping theorem
guarantees
certain natural ``periodicity''
property of the above homeomorphism
(see \eqref{periodstrip:eq} below).
This ensures that the transformed operators $H_D$ and $H_N$
have periodic coefficients.

To be precise, the quasi-conformal
Riemann's theorem
alone is not sufficient for
our needs. It is also necessary to
study the boundary behaviour
of the quasi-conformal homeomorphism.
Moreover,
we allow the boundary of
the waveguide $\Om$
to have corners and inward peaks, which give rise to
singularities of the
map.
For conformal maps this circle of questions
is exhaustively studied in the relevant literature
(see e.g. \cite{Pom}), and is usually associated with the names
of  P. Koebe, O.D. Kellogg
and S.E. Warschawski.
The extension to the quasi-conformal case
is either well known or evident to experts in
complex analysis. Nevertheless, we have
decided to provide precise
statements and complete proofs,
since we were unable to find them in the
literature in the form readily
suitable for our purposes.
These results along with the quasi-conformal
analogue of Riemann's mapping theorem
are collected in Sect. \ref{conf7:sect}.

As was mentioned earlier, the periodic isothermal coordinates
were used in \cite{KuchLev},
\cite{BSS} to establish the absolute continuity
of the operator $H$ acting in $\plainL2(\R^2)$, which
we denote below by $H_F$.
They were \underline{constructed} in \cite{KuchLev}
in the following way.
Using results from \cite{BJS}, \cite{Ve}
one can define an analytic
structure on the two-dimensional torus
with the help of local $q$-quasi-conformal coordinates.
Then it follows from the
theory of Riemann surfaces that
integration of the
analytic differential on
the torus
leads to the sought isothermal
coordinates in $\Rtwo$.
Note that the above analytic differential exists and
is unique up to a constant factor,
since the torus is a surface of genus one.
The secondary aim of our paper is to
give another, more direct proof of
the existence of such coordinates.
Instead of geometrical considerations
of \cite{KuchLev} we rely on the
well-known fact (see \cite{Ah}, \cite{AB}, \cite{BJS}, \cite{Ve})
that there exists a  \underline{unique}
$q$-quasi-conformal
homeomorphism  of $\mathbb C$ which
preserves the points $0$, $2\pi$ and  $\infty$.
Similarly to the waveguide case, a crucial observation here is that
the uniqueness combined with periodicity
of $\BG$ \underline{automatically}
implies the required periodicity of
the homeomorphism at hand
(see Theorem \ref{periodicity:thm}).
Using the isothermal coordinates, as in \cite{KuchLev}
and \cite{BSS}, we obtain the absolute continuity
of the operator $H_F$. Moreover, applying a stronger
regularity result for the coordinate change,
we are able to relax the smoothness restrictions
on $\BG$ in comparison with \cite{KuchLev}, \cite{BSS}:
for instance, under the assumption $\det\BG = 1$
the matrix $\BG$ does not need to be
H\"older-continuous, but only bounded.

After reduction to a conformal metric,
the absolute continuity of $H_F$ results immediately
from the earlier papers \cite{BS2}, \cite{BS}.
As to the operators $H_D, H_N$,
the isothermal change of variables
reduces them to operators of the same type
acting on a straight strip in $\R^2$,
with a scalar constant $\BG$.
From this point on we follow the
strategy suggested in \cite{SobWal}.
It consists in further reduction
to an auxiliary operator $H_P$ with \textsl{periodic}
conditions on the boundary of the strip of the double width.
Then a reference to \cite{BSS}
secures the required absolute continuity.

The paper is organised as follows.
Section \ref{conf2:sect} contains some preliminary material
and the precise
statements of all main results of the paper (Theorems
\ref{rtwo:thm}, \ref{dirneum1:thm}).
In Section \ref{conf3:sect} we make first reductions
simplifying the problem. In particular it is shown that
it suffices to prove the main results for the case
$\det\BG= 1$.
This section includes also
the Floquet decomposition
of the auxiliary periodic operator $H_P$.
The important Section \ref{conf4:sect} is devoted to
a detailed description of the isothermal change of variables.
Two central theorems of this section
(Theorems \ref{change:thm} and \ref{changeom:thm})
are proved in Sect. \ref{conf6:sect} after
having been translated into
the language of the quasi-conformal maps.
The proof of the main results is completed
in Sect. \ref{conf5:sect}.
The necessary information on the
quasi-conformal maps and
their boundary properties
is collected in Sect. \ref{conf7:sect}.

\section{Main results}
\label{conf2:sect}

\subsection{Notation}
\textsl{Lattices and domains.}
Let $\be_1, \be_2$
be the canonical basis in $\Rtwo$.
Along with the standard two-dimensional square lattice
$\BOLG = (2\pi\Z)^2$, introduce two ``one-dimensional''
lattices:
\begin{gather*}
\bg_1 = (2\pi\Z)\times\{0\}
= \{ 2\pi n\be_1,  n\in\Z\},\\
\bg_2 = \{0\}\times (2\pi \Z)
= \{ 2\pi n \be_2, n\in\Z \}.
\end{gather*}
We say that a function $f$ is $\bg_j$-periodic (resp.
$\BOLG$-periodic),
if $f(\bx+ 2\pi n \be_j) = f(\bx)$ a.a. $\bx$ and all $n\in\Z$
(resp. $f(\bx+ \bxi) = f(\bx)$ a.a. $\bx$ and all
$\bxi\in\BOLG$).

For any set $\CF\subset\Rtwo$ define its translates
as follows:
\begin{gather*}
\CF^{(\bn)} = \{\bx\in\Rtwo: \bx-2\pi\bn\in\CF\},\ \bn\in\Z^2,\\
\CF^{(n)} = \CF^{(\bn)},\ \textup{with}\ \ \bn = (n, 0),\ n\in\Z.
\end{gather*}
We say that the set $\CF$ is $\BOLG$-\textsl{periodic}
if $\CF = \CF^{(\bn)}$ for all $\bn\in \Z^2$,
and that $\CF$ is $\bg_1$-\textsl{periodic} if
$\CF = \CF^{(n)}$ for all $n\in \Z$. When it does not lead
to a confusion, instead of ``$\bg_1$-'' or ``$\BOLG$-periodicity''
we use the term ``periodicity''.
Similarly we define the periodicity of sets $\CF$
on the cylinder
\begin{equation*}
\CC = \Rtwo/ \bg_2.
\end{equation*}
Precisely, $\CF\subset\CC$ is said to be periodic
(or $\bg_1$-periodic) if $\CF^{(n)} = \CF$ for all $n\in\Z$.

We are going to study periodic operators
in two-dimensional domains of three types:
on the entire plain $\Rtwo$,
on the cylinder $\CC$ or on a
periodic domain $\Om\subset\Rtwo$.
Our model periodic domain will be
the straight  strip
\begin{equation}\label{strip:eq}
\CS = \CS_{\dc} =
\{\bx = (x_1, x_2)\in \Rtwo : 0< x_2 < \pi \dc\}, \ \dc >0.
\end{equation}
Obviously, the cylinder $\CC$
can be viewed as the closed strip
$\overline{\CS_2}$ with identified lower and upper boundaries.
In general, we assume that $\Om$ is as
described below:

\begin{defn}\label{periodic:defn}
We say that a domain $\Om\subset\Rtwo$
is admissible if
there exists a finite collection of
bounded domains $\CE_j,\ j = 1, 2, \dots, N$
with Lipschitz boundaries
such that the set
\begin{equation}\label{ce0:eq}
\CE_0 = \bigcup_{j = 1}^N \CE_j
\end{equation}
is connected, and $\Om = \cup_{n\in\Z} \CE_0^{(n)}$.
\end{defn}

Note that the domain $\Om$ is automatically
$\bg_1$-periodic and bounded in the direction $\be_2$.
We also point out that $\Om$ satisfies
the cone condition (see e.g.
\cite{MP}, Sect. 1.3.3 ),
since so do the domains $\CE_j, j =1, 2, \dots, N$.
Without loss of generality we \textbf{always assume that}
\begin{equation}\label{found:eq}
\CE_0\supset
\{\bx\in \Om: -1\le x_1 \le 2\pi +1 \}.
\end{equation}
Certainly, the choice of the domains $\CE_j$ for a given
admissible $\Om$ is not unique.

As a rule we  try to treat all three
cases simultaneously,
and therefore we use the notation $\L$ either for
$\Rtwo$ or $\CC$ or $\Om$.
Identifying the points that differ by a vector of
the lattice we define
\begin{equation}\label{u:eq}
\U =
\begin{cases}
\Om/\bg_1,\ \ \textup{if}\ \ \L = \Om;\\[0.2cm]
{\Ttwo = \Rtwo/{\BOLG}},\ \ \textup{if}\ \
\L = \Rtwo \ \ \textup{or} \ \ \CC.
\end{cases}
\end{equation}
Introduce also the \textsl{fundamental domains}:
\begin{equation*}
\CO =
\begin{cases}
{(0, 2\pi)\times (0, 2\pi)},\ \ \textup{if}\ \
\L =  \Rtwo \ \ \textup{or}\ \ \ \CC,\\
\{\bx\in\Om : 0 < x_1 < 2\pi\},\ \ \
\textup{if} \ \ \ \L = \Om.
\end{cases}
\end{equation*}
In the case $\L = \Om$ the set $\CO$
might not be connected.

Similarly to $\Om$, for
our purposes it will be necessary
to view $\R^2$ as being covered by bounded domains.
Precisely, we assume that there are finitely many
bounded domains $\CE_j, j = 1, 2, \dots, N$,
with Lipschitz boundaries such that
\begin{equation*}
\overline{\CO}\subset \CE_0 = \bigcup_{j=1}^N \CE_j.
\end{equation*}

\textsl{System of curves.}
We also need to introduce a system of
curves $\BSIG$ in $\L$ associated with the
covering of $\L$ by $\CE_j$'s and their translates.
Below by a ``$\plainC{m+\a}$-arc", $m\in \mathbb N$,
\ $0\le \a <1 $,
we mean a Jordan arc in $\Rtwo$ which is
parametrised by a $\boldC{m+\a}$-smooth function
$\boldsymbol\psi:[0, 1]\to \Rtwo$ such that
$|\boldsymbol\psi'(t)|> 0,\ t\in[0, 1]$.

\begin{defn}\label{sigma:defn}
Let $\L = \Rtwo$ or $\Om$, and let the parameter
$\bn$ vary over the set $\Z^2$ ( for $\L = \R^2$)
or $\Z\times\{0\}$ (for $\L = \Om$).
Let $\ell_{j}\subset \overline{\CE_j},\ j = 1, 2, \dots N$,
be a finite set
of $\plainC1$-arcs such that
\begin{equation*}
\ell_{j}^{(\bn)}\cap \ell_{j} = \emptyset,\
\forall \bn\not = 0.
\end{equation*}
Then \textsl{a (periodic)
system of curves} $\BSIG$ \textsl{in} $\L$
is defined to be a \textsl{family} $\{\Sigma_{j}\}$
of the closed sets $\Sigma_{j}\subset \overline{\L}$
of the form
\begin{equation}\label{sigmaj:eq}
\Sigma_{j} = \cup_{\bn} \ell_{j}^{(\bn)}.
\end{equation}

Let $\L = \CC$. Then the set $\BSIG$ is
\textsl{a system
of curves in} $\CC$ if there exists a system of curves
$\boldsymbol\Mu$ in $\Rtwo$ such that $
\BSIG = \boldsymbol\Mu/\bg_2$, where the latter being defined as
the family $\{\Mu_{j}/\bg_2\}$.

Similarly one defines
\begin{equation}\label{xi:eq}
\BXI =
\begin{cases}
\BSIG/\bg_1,\ \ \textup{if}\ \
\L = \Om\ \ \textup{or} \ \ \CC;\\[0.2cm]
\BSIG/{\BOLG},\ \ \textup{if}\ \
\L = \Rtwo,\
\end{cases}
j = 1, 2, \dots, N.
\end{equation}
\end{defn}

\vskip 0.3cm

\begin{rem}\label{curves:rem}
\begin{itemize}
\item[(i)]
In the above definition any two curves
$\ell_{j}$, $\ell_{l}$
or their parts can coincide.
This means in particular that $\BSIG$
may contain several ``copies'' of the same curve.
The curves  are also allowed to meet at the zero angle.

\item[(ii)]
Definition \ref{sigma:defn}
prescribes exactly one curve $\ell_j$ for
each of the domains $\CE_j$. A seemingly more general
situation when there are finitely many curves in
$\overline{\CE_j}$ is easily reduced to the original one
by adding a necessary number of copies of $\CE_j$
in \eqref{ce0:eq}.
\end{itemize}
\end{rem}

For any bounded set
$\CE\subset \overline{\L}$ we denote
\begin{equation}\label{sigmace:eq}
\BSIG_{\CE} = \{\Sigma_{j,\CE}\},\ \
\Sigma_{j, \CE} = \cup_{\bn} \ell_{j}^{(\bn)},
\end{equation}
where $\bn$ are such that
$\ell_j^{(\bn)}\cap\CE\not=\emptyset$.
Obviously, each $\Sigma_{j, \CE}$
contains finitely many translates
of $\ell_{j}$.

\textsl{Function spaces.}
Apart from the standard
notational conventions for classes of
differentiable functions we use
the notation $\plainW s p (\CC)$, $s = 0, 1$,\ $p\ge 1$
which stands for the space of all
$\bg_2$-periodic functions from
$\plainW s p_{\loc}(\Rtwo)$ equipped with the norm
$\|\ \cdot\ \|_{\plainW s p(\CS_2)}$.
Similarly, the notation $\plainW s p(\U)$,
$s = 0,1$,\  $p\ge 1$,  means  the
space of all $\bg_1$-periodic ( for $\L = \Om$)
or $\BOLG$-periodic (for $\L = \Rtwo$ or $\L = \CC$)
functions $u\in\plainW s p_{\loc}(\L)$
with the norm $\|\ \cdot\ \|_{\plainW s p(\CO)}$.
In the case $p = 2$ we use the standard notation
$\plainH s \equiv \plainW s 2$.
The same convention applies to H\"older spaces
$\plainC{m+\a}$.
For instance
$\plainC{m+\a}(\U)$, $m\in \mathbb N\cup \{0\},\ \a\in (0, 1)$
(resp.
$\plainC{m+\a} (\overline{\U})$) denotes the space of
all $\bg_1$-periodic ( for $\L = \Om$)
or $\BOLG$-periodic (for $\L = \Rtwo$ or $\L = \CC$)
functions $u\in\plainC{m+\a}(\L)$
(resp. $\plainC{m+\a}(\overline{\L})$ ). Certainly,
in the case $\L = \Rtwo$ or $\L = \CC$ the space
$\plainC{m+\a}(\U)$ coincides with
$\plainC{m+\a}(\overline{\U})$.

To denote spaces of vector-valued functions, we use
the boldface letters, e.g. $\boldL p(\L)$.
The notation $\GD\boldf$ stands for the Jacobian
matrix of the function $\boldf$.
Slightly abusing the notation we sometimes do
not distinguish between $\plainL2(\CC)$ and $\plainL2(\CS_2)$
($\plainL2(\U)$ and  $\plainL2(\CO)$).

Function spaces on $\BSIG$ are defined in a natural way.
Namely, by definition the space $\boldL p(\BSIG)$
is a set of functions $\bsig = \{\s_{j}\}$
such that
$\s_{j}\in \plainL p(\Sigma_{j}), j = 1, 2, \dots, N$.
Similarly one defines  $\boldL p(\BXI)$.
We say that $\bsig$ is real-valued if
all the components $\s_{j}$ are real-valued.

Definition of traces on $\BSIG$
of functions on $\L$ requires
special comments. Suppose that $f\in \plainH1(\L)$.
The trace $f|_{\ell_{j}}$ is defined to be the trace
of the function $f|_{\CE_j}\in \plainH1(\CE_j)$.
Sometimes we write $f$
instead of the trace $f|_{\ell_j}$ when it
does not cause confusion.
By the embedding theorems and
multiplicative inequality
for the traces
(see \cite{Maz}, Corollary 1.4.7/2),
$f\in \plainL r(\ell_j)$
with any $r <\infty$ and
\begin{equation}\label{traces:eq}
\|f\|_{\plainL r(\ell_j)}
\le \e \|\nabla u\|_{\plainL2(\CE_j)}
+ C_j(\e) \|u\|_{\plainL 2(\CE_j)},
\end{equation}
for all $\e >0$.
Similarly one defines the traces on the translates
$\ell^{(\bn)}_{j}$, which leads to the
collection of traces
$f|_{\BSIG} = \{f|_{\Sigma_{j}}\}$
in a natural way.
Note that if an arc $\ell = \ell_j = \ell_m$
belongs to $\CE_j\cap \CE_m$,
then $f|_{\ell_j} = f|_{\ell_m}$.
On the contrary,
if $\ell\subset \p\CE_j\cap\p\CE_m$ for two
distinct $j$ and $m$, then $f|_{\ell_j}$ may be different
from $f|_{\ell_m}$.
The latter situation can occur
in the case
when $\L = \Om$ and $\ell$ is a part of the boundary of $\Om$
such that $\Om$ lies on both sides of $\ell$.

\textsl{Weights and coefficients.}
We work in the weighted space $\plainL2(\mu, \L)$ with the norm
\begin{equation*}
\|u\|_\mu = \biggl[\int_\L |u|^2 \mu d\bx\biggr]^{1/2},
\end{equation*}
where $\mu$ is a real-valued periodic
function, satisfying the conditions
\begin{equation}\label{weight:eq}
\begin{cases}
\mes\{ \bx\in\L: \mu(\bx)\le 0\} = 0,\\
\mu\in \plainL t(\U),\
t >1.
\end{cases}
\end{equation}
In the cases when it does not cause any
confusion, the subscript $\mu$ will be
omitted from the notation of the norm.
We are interested in
the properties of
the Schr\"odinger type
periodic operators in $\plainL2(\mu, \L)$,
defined by the formal expression
\begin{equation*}
H  = \frac{1}{\mu}
\lu (\BD-\ba),\om^2\BG(\BD-\ba)\ru + \frac{1}{\mu}V
+ \bsig\d_{\BSIG}, \
\BD = -i\nabla,
\end{equation*}
where $V, \om, \ba, \BG$ are some real-valued periodic
(vector/matrix) functions defined on $\L$, and
$\bsig$ is a real-valued function
defined  on a system of
curves $\BSIG$.
Let us now give a precise definition of the
operators in question.
The \textsl{electric potential}  $V$,
the \textsl{magnetic vector-potential} $\ba$
and the function $\bsig$
are supposed to satisfy the conditions
\begin{gather}
\ba\in \boldL{s}(\U),\
 s > 2;\ \ V\in \plainL{p}(\U),
\ p >1,\label{ndcond:eq}\\[0.3cm]
\bsig\in\boldL{r}(\BXI),\  r >1.\label{s:eq}
\end{gather}
The coefficient $\BG = \{g_{jl}(\bx)\}, j,l = 1, 2$
is a symmetric matrix-valued function on $\U$
with real-valued entries
$g_{jl}(\bx)$ that satisfy the condition
\begin{equation}\label{g:eq}
\begin{cases}
c|\bxi|^2\le \langle\BG(\bx)\bxi,
\bxi\rangle\le C|\bxi|^2,
\\[0.2cm]
\det\BG(\bx) = C',
\end{cases}
\ \ \ \forall \bxi \in
\Rtwo,\ \textup{a.a.} \ \   \bx\in \U.
\end{equation}
Here and everywhere below
by $C$ and $c$ with or
without indices we denote various positive
constants whose precise value is unimportant.
As a rule we assume that $\det\BG = 1$, but
sometimes it is convenient not to have this restrictions.

As to $\om$, it is a real-valued function on $\U$, such that
\begin{equation}\label{ommu:eq}
c\le \om(\bx)\le C, \ \ \textup{a.a.}\ \  \bx\in \U.
\end{equation}
We are interested in four different realisations of the
operator $H$. Namely, we study
\begin{itemize}
\item
$H$ as an
operator in $\plainL2(\mu, \Rtwo)$,
which will be later referred to
as the ``full'' operator $H_F$;
\item
$H$ as an operator in $\plainL2(\mu, \CC)$,
or which is the same, as an operator
in $\plainL2(\mu, \CS_2)$
with periodic conditions on the boundary of $\CS_2$.
In this case we use the notation
$H_P$;
\item $H$ acting in $\plainL2(\mu, \Om)$ with the
Dirichlet or natural boundary condition.
In these cases we use the notation $H_D$ or $H_N$ respectively.
\end{itemize}
When we need to treat all four situations
simultaneously, we use the notation
$H_\aleph$, where $\aleph$ will have the meaning of any of the
four letters $F, P, D, N$.

For each of these four problems
the operator will be defined via its
quadratic form.
To give this definition it suffices to
assume \eqref{weight:eq}, \eqref{ndcond:eq}, \eqref{s:eq},
\eqref{g:eq} and \eqref{ommu:eq},
although later on we shall need more restrictive
conditions on $\BG$ and $\om$.
Consider the following quadratic form
\begin{equation}\label{qform:eq}
h[u] =
\int_{\L} \om^2
\langle\BG(\BD - \ba)u, \overline{(\BD - \ba) u}\rangle
d\bx
+ \int_{\L} V |u|^2 d\bx
+ \int_{\BSIG} \bsig |u|^2 dS,
\end{equation}
where
\begin{equation}\label{sig:eq}
\int_{\BSIG}\bsig|u|^2 dS =
\sum_{j=1}^N \int_{\Sigma_j}\s_j \bigl|u|_{\Sigma_j}\bigr|^2 dS
\end{equation}
defined  either on
$\CD_F = \plainH1(\Rtwo)$
(for the full problem), or
 $\CD_P = \plainH1(\CC)$
(for the $\bg_2$-periodic problem),
or
$\CD_D = \plainH1_0(\Om)$ or $\CD_N = \plainH1(\Om)$.
Depending on the domain
we denote the form $h$ by
$h_F, h_P, h_D$ or $h_N$ respectively.
Sometimes, in order to distinguish operators
defined on different domains
or/and with different coefficients and weights,
we use the full notation
$H_{\aleph}(\om, \BG, \ba,  V, \mu, \bsig; \L)$
or its short
variants: e.g. $H_{\aleph}(\om, \BG)$.
Similar convention applies to the notation of the
quadratic forms $h_{\aleph}$.

Let us check that these forms are closed
in $\plainL2(\mu, \L)$. To this end
split the form \eqref{qform:eq}
into the \textsl{unperturbed} form and the
\textsl{perturbation} form:
\begin{align*}
h_\aleph[u] = &\ h^{(0)}_{\aleph}[u]
+ w_{\aleph}[u];\\
h^{(0)}_{\aleph}[u] = &\ \int_{\L} \om^2\lu\BG
\BD u, \overline{\BD u}\ru d\bx, \\
w_{\aleph}[u] = &\ \int_{\L}\bigl[
- \om^2\lu \BG \BD u, \ba \ru \overline{u}
-  \om^2 u\lu \BG \ba, \overline{\BD u} \ru
+ \om^2\lu\BG\ba, \ba\ru
|u|^2\bigr] d\bx\\
&\qquad\qquad + \int_{\L} V |u|^2 d\bx
+ \int_{\BSIG} \bsig |u|^2 dS.
\end{align*}
The necessary properties of these forms are
contained in

\begin{lem}\label{qforms:lem}
Let $\mu$ satisfy the condition \eqref{weight:eq} and let $\ba$,
$V$ and $\bsig$ satisfy \eqref{ndcond:eq}, \eqref{s:eq}. Then
\begin{itemize}
\item[(i)]
The standard $\plainH1$-norm is equivalent to the norm induced by
the form $h^{(0)}_{\aleph}$, i.e.
\begin{equation}\label{qforms:eq}
C^{-1}(\|\nabla u\|_1^2 + \|u\|_1^2)
\le h^{(0)}_{\aleph}[u]
+ \|u\|_\mu^2
\le C(\|\nabla u\|_1^2 + \|u\|_1^2),\ \forall u\in \CD_\aleph,
\end{equation}
with some constant $C$ depending on $\om, \BG, \mu$;
\item[(ii)]
The forms $h^{(0)}_\aleph$ with the domains $\CD_\aleph$ are closed in
$\plainL2(\mu, \L)$;
\item[(iii)]
For any $\e>0$ there exists a constant $C_\e$ such that
\begin{equation}\label{infinite:eq}
|w_{\aleph}[u]|\le \e h^{(0)}_{\aleph}[u] + C_\e \|u\|_\mu^2,\
\forall u\in \CD_\aleph,
\end{equation}
so that the perturbed forms
$h_\aleph$ are also closed.
\end{itemize}
\end{lem}

\begin{proof} By virtue of
\eqref{g:eq} and \eqref{ommu:eq}, we may assume without loss of
generality that $\BG = \BI$ and $\om = 1$.

To be definite, we shall establish the equivalence
\eqref{qforms:eq} for the case $\L = \Om$ only.
Before doing this, we remind the following general fact.
Let $\CE\subset\Om$ be an arbitrary bounded domain
satisfying the cone
property. By \eqref{weight:eq}
the equality
\begin{equation*}
\int_\CE |u|^2 \mu d\bx = 0
\end{equation*}
for a constant function $u$ implies
that $u = 0$ in $\CE$.
Since $\plainH1(\CE)$ is embedded into
$\plainL{h}(\CE)$ for any $h<\infty$
(see e.g. \cite{MP}, Sect. 1.8.2),
we see that the functional in
the l.h.s. of the last inequality
is bounded by
\begin{equation*}
\|\mu\|_{\plainL t(\CE)} \| u\|^2_{\plainL{h}(\CE)}
\le C_\mu \|u\|_{\plainH1(\CE)}^2,
\ \ \frac{2}{h} + \frac{1}{t} = 1,
\end{equation*}
i.e. it is continuous in $\plainH1(\CE)$.
Consequently, (see e.g. \cite{MP}, Sect. 1.5.4) the functional
\begin{equation*}
\int_\CE |\nabla u|^2 d\bx + \int_{\CE} |u|^2 \mu d\bx
\end{equation*}
defines in $\plainH1(\CE)$ a norm equivalent to the
standard $\plainH1(\CE)$-norm.

Let $\CE_0$ be the domain defined in \eqref{ce0:eq}.
By Definition \ref{periodic:defn}
the translates $\CE_0^{(n)}, n\in \Z$  form
an open covering of $\Om$, and each of these domains
satisfies the cone condition.
Let us construct a partition
of unity subordinate to the covering of $\Om$ by
the sets $\CE_0^{(n)},\ n\in\Z$.
Let $\z\in \plainC\infty_0(\R)$ be a function
such that  $0\le \z(t) \le 1$,
$\z(t) = 0,\ t\notin (-1, 2\pi +1)$, and
\begin{equation}\label{partition:eq}
\sum_{n\in\Z}\z_n^2 = 1,\ \z_n(t) = \z(t-2\pi n).
\end{equation}
From \eqref{found:eq} it is clear that for any $u\in\plainH1(\Om)$
the truncated functions
$u_n(\bx) = u(\bx)\z_n(x_1)$ belong to
$\plainH1(\CE_0^{(n)})$,\ so that all $u_n$'s satisfy
the equivalence relation of the form
\eqref{qforms:eq}, with a constant
$C$ independent of $n$.
Summing them over $n\in\Z$
we obtain from \eqref{partition:eq}
and the equality $\sum_n \z_n\nabla\z_n = 0$ that
\begin{equation*}
\sum_n \biggl[
\|\nabla u_n\|^2_1
 + \|u_n\|^2_\mu \biggr]
= \|\nabla u\|^2_1 + \|u\|^2_\mu
+ \int_{\Om}\sum_{n}(\p_1\z_n)^2 |u|^2 d\bx.
\end{equation*}
Since the sum in the last
integral is bounded,
this leads to the required equivalence in $\Om$.

To prove (ii)
it remains to show that $\CD_\aleph$ is embedded in
$\plainL2(\mu, \L)$. By the first part of the proof
any function $u\in \CD_\aleph\subset \plainH1(\L)$
belongs to $\plainL2(\mu, \L)$.
By virtue of \eqref{weight:eq}
if $u=0$ $\mu$-a.e., then
inevitably $u = 0$ a.e. with respect to
the standard Lebesgue measure.

To prove (iii) we use the same partition of unity as in the
first part of the proof. Since $\sum_n \z_n \nabla\z_n = 0$,
we have
\begin{equation*}
\sum_n w_\aleph[u_n] = w_\aleph[u],
\end{equation*}
and hence the problem reduces to
checking the estimates \eqref{infinite:eq}
for each $u_n$ individually.
Furthermore,
in view of the periodicity of the domain, it suffices
to prove \eqref{infinite:eq}
for the domain $\CE_0 = \CE_0^{(0)}$
only. Below we denote $u = u_0$.

Suppose first that $\bsig = 0$.
Let us fix an $\e >0$ and define
\begin{equation*}
\CE'_\d = \{\bx\in\CE_0: \textup{either}\ \
|\ba(\bx)|>\d^{-1}, \ \ \textup{or}\ \
|V(\bx)|>\d^{-1}
\ \textup{or}\ \ \mu(\bx)<\d\}.
\end{equation*}
Denote by $\chi_\d$ the characteristic function
of $\CE'_\d$.
Since the functions $V_1 = (1-\chi_\d)V,\ \ba_1 = (1-\chi_\d)\ba$
are bounded and $\mu(1-\chi_\d)\ge \d(1-\chi_\d)$,
it is straightforward to check, using H\"older's
inequality, that
\begin{align*}
|w_\aleph(\ba_1, V_1, \mu; \CE_0)[u]|
\le &\ \e\|\nabla u\|_1^2
+ C_\d(4\e)^{-1}\|u(1-\chi_\d)\|_1^2\\
\le &\ \e\|\nabla u\|_1^2
+ C_\d(4\e\d)^{-1}\|u\|_\mu^2,
\end{align*}
for any $\e >0$.
On the other hand, using again
the embedding of $\plainH1(\CE)$ into
$\plainL h(\CE)$ for any $h<\infty$ (see \cite{MP}, Sect. 1.8.2),
we easily obtain that
\begin{equation}\label{chi:eq}
|w_\aleph(\ba\chi_\d, V\chi_\d, \mu; \CE_0)[u]|
\le C\bigl(\| \ba\chi_\d\|_{\plainL s}
+ \| \ba\chi_\d\|_{\plainL s}^2 + \| V\chi_\d\|_{\plainL p}
\bigr)
(\|\nabla u\|_1^2 + \|u\|_1^2).
\end{equation}
The norms of $\ba\chi_\d, V\chi_\d$
can be made arbitrarily small
by choosing $\d$ appropriately.
Now the r.h.s. of \eqref{chi:eq} satisfies the bound
\eqref{infinite:eq} by virtue of \eqref{qforms:eq}.

It remains to establish the bound \eqref{infinite:eq}
for the form $w_\aleph$ with $\ba = 0$, $V = 0$ and
an arbitrary $\bsig$, i.e. for
\begin{equation*}
w_\aleph[u] = \int_{\BSIG} \bsig |u|^2 dS
= \int_{\BSIG_{\CE_0}}\bsig |u|^2 dS,\ \  u = u_0,
\end{equation*}
where the set $\BSIG_{\CE_0}$ is defined by
\eqref{sigmace:eq}. Since this set contains only finitely
many translates of each arc $\ell_j$, it suffices to establish the
sought estimate for each of them individually.
Denote $\ell = \ell^{(n)}_j$, $\s  = \s_j$
with some $j = 1, 2, \dots, N$
and $n\in \Z$. Then by H\"older's
inequality and \eqref{traces:eq},
for any $\e>0$ we have
\begin{equation*}
\int_{\ell} \s |u|^2 dS
\le \|\s\|_{\plainL r(\ell)} \ \|u\|_{\plainL t(\ell)}^2
\le \e \int_{\CE_j} |\nabla u|^2 d\bx + C(\e, \s)\int_{\CE_j}
|u|^2 d\bx.
\end{equation*}
By \eqref{infinite:eq} with $\bsig = 0$ the
last term is bounded from above by
\begin{equation*}
\e\|\nabla u\|_1^2 + \tilde C(\e, \s)\|u\|_{\mu}^2.
\end{equation*}
The required result follows.
\end{proof}

By Lemma \ref{qforms:lem}
all four forms $h_\aleph$ are closed, and therefore
they uniquely define in $\plainL2(\mu, \L)$
four self-adjoint operators which we denote
by $H_F$, $H_P, H_D$ and $H_N$.
We do not need to know the domains of these operators,
although they can be specified under
supplementary regularity
conditions on the coefficients,
the boundary of $\Om$ and the
curves from $\BSIG$.
For instance, if the system of curves
contains only the boundary $\p\Om$, then
$H_N$ is the operator of the third boundary value problem with the
condition
\begin{equation*}
\om^2\lu \BG (\nabla - i\ba) u, \bn\ru + \s u = 0,\  \bx\in \p\Om,
\end{equation*}
where $\bn = \bn(\bx)$
is the exteriour unit normal to the boundary at $\bx\in\p\Om$.
If $\BSIG$ contains a curve which has an arc $\ell$ strictly inside
$\Om$ and separated from other components of $\BSIG$, then
the integral over $\BSIG$ in \eqref{qform:eq} induces
the condition
\begin{equation*}
\bigl[\om^2\lu \BG \nabla u, \bn\ru\bigr] + \s u = 0,\ \bx\in\ell,
\end{equation*}
on the jump $[\dots]$ of the conormal derivative across the
curve $\ell$.

An important role  will be played
by the general observation that
the singular continuous spectra of the operators
$H_\aleph$ are empty, which we state separately for
later reference:

\begin{prop}\label{sc:prop}
Let $\L$ be either $\Rtwo$ or $\CC$ or $\Om$.
Suppose that the conditions
\eqref{weight:eq}, \eqref{ndcond:eq},
\eqref{g:eq}, \eqref{ommu:eq} and \eqref{s:eq}  are fulfilled.
Then the singular continuous spectra of
$H_\aleph$ are empty.
\end{prop}

The proof of this property is based
on the standard direct integral
representation for $H_\aleph$, known as
\textsl{the Floquet decomposition}.
The crucial fact is that the resolvents of the
fibers of $H_\aleph$
in this representation are compact
operator-functions, analytic in
the quasi-momentum. We are not going
to provide all the details of this argument,
but refer to the comprehensive exposition of this
issue in \cite{Kuch}, and also to \cite{KuchLev},
\cite{HH} and \cite{W}.
The  Floquet decomposition for $H_F$ can be found,
for example, in \cite{BSS} and \cite{BS}.
For the operators $H_D, H_N$
it is clearly explained in \cite{BirSol}.
For our purposes
we need to describe the Floquet decomposition
for the operator $H_P$ only
(see Sect. \ref{conf3:sect}).

\subsection{Results}
Our main goal is to go further and prove that the spectrum of
$H_\aleph$ is absolutely continuous. Here $\aleph$ takes the
values $F, D$ or $N$. The operator $H_P$ plays an auxiliary role:
we prove its absolute continuity only in the case of a diagonal
constant matrix $G$. This result will be decisive in the proofs
for the cases $D$ and $N$.

Now the conditions \eqref{g:eq},
\eqref{ommu:eq} are insufficient.
Suppose in addition that
\begin{equation}\label{om:eq}
\om\in \plainW1q(\U) \ \
\textup{and}\ \ \
\BG\nabla\om\in \boldW1{q}(\U),\ q > 1.
\end{equation}
As far as $\BG$ is concerned,
a number of results will be obtained under the
additional restriction
\begin{equation}\label{holder:eq}
\BG\in \boldC\a\bigl(\overline{\U}\bigr),\
 \ \a \in(0, 1).
\end{equation}
Clearly, for a uniformly Lipshitz matrix $\BG$ the condition
\eqref{om:eq} is equivalent to $\om\in\plainW2q, q > 1$.
One is tempted to say that due to the
presence of the function $\om$ in
\eqref{qform:eq},
the condition $\det\BG = \textup{const}$
in \eqref{g:eq} does not restrict
generality.
We emphasise however that the smoothness
conditions for the functions $\BG$ and $\om$
are different.

The main result for the operator
$H_F$ is
contained in the next theorem.

\begin{thm}\label{rtwo:thm}
Let $\L = \Rtwo$, and $\bsig = 0$.
Suppose that the conditions
\eqref{weight:eq}, \eqref{ndcond:eq},
\eqref{g:eq}, \eqref{ommu:eq}
and \eqref{om:eq} are fulfilled.
Then the spectrum of $H_F$ is absolutely continuous.
\end{thm}

\begin{rem}\label{constant:rem}
For a variable $\BG$ satisfying
\eqref{holder:eq} Theorem
\ref{rtwo:thm} was proved in
\cite{BSS} even with a
$\bsig\not = 0$.
R. Shterenberg \cite{Sh2}
recently announced the absolute continuity
of $H_F$ under the reduced smoothness assumption
$\om\in\plainW 1 q$, $q > 2$.
For our proof  we need an earlier
result from the paper \cite{BS2} where
Theorem \ref{rtwo:thm} was established for
constant matrices $\BG$ and $\om = 1$.
\end{rem}

To prove the absolute continuity of the
operators $H_D$, $H_N$
we need to impose some extra conditions on
the domain $\Om$.
We assume that $\Om$ is $\bg_1$-periodic and that
there exists a homeomorphism of
$\overline{\Om}$ onto $\overline{\CS_1}$. This means,
in particular, that the boundary of $\Om$ consists
of two disjoint $\bg_1$-periodic
Jordan curves that we denote by $\ell_-, \ell_+$.
We impose the following condition on these curves:

\begin{cond}\label{boundary:cond}
\begin{enumerate}
\item
Locally, each curve $\ell_-, \ell_+$,
is a piece-wise $\plainC{1+\a}$-smooth Jordan arc
with $\a\in (0, 1)$;
\item
The domain $\Om$ does not have any outward peaks,
i.e. the interior angle between any two smooth
components of the boundary
at each point of non-smoothness is strictly greater
than zero.
\end{enumerate}
\end{cond}

Note that Condition \ref{boundary:cond} does not exclude
\textsl{inward} peaks, which guarantees that
that the cone condition is satisfied. Using this fact
and remembering that every bounded domain with cone condition
can be represented as a union of finitely many domains
with Lipshitz boundaries (see e.g. \cite{Maz}),
one can easily show that the domain $\Om$ is admissible
in the sense of Definition \ref{periodic:defn}.

The parameter $\a$
above, without loss of generality can
be chosen the same as in
the condition \eqref{holder:eq}.
Below we denote by $\bn(\bx)$, a.a. $\bx\in\p\Om$,
the exteriour unit normal to the
boundary $\p\Om$, and by $Z\subset\p\Om$
the discrete set  where the smoothness of the boundary
breaks down.
Clearly, $Z$ has no finite accumulation points.

\begin{thm}\label{dirneum1:thm}
Let $\L = \Om$.
Suppose that the domain $\Om$ satisfies
Condition \ref{boundary:cond}, the set
$\BSIG$ is a system of curves on $\Om$
in the sense of Definition \ref{sigma:defn},
and the conditions
\eqref{weight:eq}, \eqref{ndcond:eq},
\eqref{g:eq},  \eqref{ommu:eq},  \eqref{s:eq},
\eqref{om:eq}, \eqref{holder:eq}
are fulfilled.
Then the spectra of $H_D$ and $H_N$
are absolutely continuous.
\end{thm}

For later convenience we make a couple of
simplifying assumptions that do not
restrict generality.

First we include the
boundary $\p\Om$ in the
system of curves $\BSIG$,
even if $\BSIG$ already contains either pieces of $\p\Om$
or the entire boundary.
The notation for the components
of $\BSIG$ will be as follows.
By Condition \ref{boundary:cond} there are finitely many
$\plainC{1+\a}$-arcs $\ell_j\subset\p\Om$,
\ $j = 1, 2, \dots, M<\infty$
such that
\begin{equation*}
\p\Om = \cup_{j=1}^M\Sigma_j,
\end{equation*}
where $\Sigma_j$ are defined by \eqref{sigmaj:eq}.
Since the boundary $\p\Om$ is a Jordan curve,
we can assume that any pair of arcs $\ell_j, \ell_s$
with $s,j = 1, 2, \dots, M$ do not have
common interiour points.
Moreover, since the outward peaks are absent, there are
$M$ bounded domains $\CE_j\subset\Om$ with Lipschitz
boundaries such  that $\ell_j\subset\overline{\CE_j}$.
This allows us to include $\Sigma_j, j = 1, 2, \dots, M$
in the initial system of curves in $\Om$.
From now on we consider the sets
$\Sigma_j, j = 1, 2, \dots, M$ to
be the first $M$ components of $\BSIG$, and the
original components of $\BSIG$
will be relabelled to have numbers
from $M+1$ to $N$.
Also we set $\s_j = 0$,\ $j = 1, 2, \dots, M$. Obviously,
this procedure does not change
the operators $H_D, H_N$, since
the quadratic form
\eqref{qform:eq} remains unchanged.

Secondly, we assume that if a
curve $\ell_j, \ j = 1, 2, \dots, N$, contains
a point $\bz\in Z$, then $\bz$ is either the start
or the end point of $\ell_j$, i.e. given a parametrisation
$\boldsymbol\psi_j: [0, 1]\to \ell_j$, we have
$\boldsymbol\psi_j(0) = \bz$ or $\boldsymbol\psi_j (1) = \bz$.
This can be done by breaking, if necessary,
every $\ell_j$, containing a $\bz\in Z$, into subarcs, and
using Remark \ref{curves:rem}(ii).

\section{Preliminary conclusions}
\label{conf3:sect}

\subsection{Reduction to $\mu=\om=1$}
Proposition \ref{sc:prop}
allows one to show that it suffices
to prove Theorems \ref{rtwo:thm},
\ref{dirneum1:thm} for $\om = \mu = 1$. The
following Lemma is a variant of
a well-known result (see e.g. \cite{BS}), and it is
a crucial ingredient in our argument.

\begin{lem} \label{bs:lem}
In addition to the conditions of
Proposition \ref{sc:prop}
assume also that \eqref{om:eq} is satisfied.
In the case of the operator $H_N$ assume also that
Condition \ref{boundary:cond} is fulfilled.
Then
\begin{gather}
\om^{-1}
H_{\aleph}(\om, V, \bsig)
\om^{-1} =
H_{\aleph}(1, \tilde V, \tilde{\bsig}),\label{sandw:eq}\\
\tilde V = \om^{-2}V
+ \om^{-1}\lu\nabla, \BG\nabla\ru \om,\notag\\
\tilde\s_j =
\begin{cases}
 -\om^{-1}\lu\BG\nabla\om, \bn\ru,\ \ \textup{if} \ \
j = 1, 2, \dots, M, \notag\\
\om^{-2}\s_j,\ \ \ \ \textup{if} \ \ j = M+1, \dots, N. \notag
\end{cases}
\end{gather}
\end{lem}

\begin{proof} We prove the Lemma for the case $\aleph = N$ only.
To avoid cumbersome calculations assume that
$\ba = 0, V = 0, \bsig = 0$.
The general case requires only obvious modifications.
The second condition in \eqref{om:eq} and \eqref{g:eq}
imply that $\nabla\om\in\plainL{h}(\U)$ with some $h>2$.
Therefore $\om\in\plainW1h(\U)\subset \plainC{}(\U)$, and
using  \eqref{ommu:eq}, one can show that
the functions
$\om$, $\om^{-1}$
are multipliers in $\plainH1(\Om)$.
This implies that
the quadratic forms of the operators in the r.h.s. and l.h.s.
of \eqref{sandw:eq}
are both closed on $\plainH1(\Om)$.
Thus it suffices to prove that the
corresponding bilinear forms coincide.
Let us consider the
form of the operator
in the l.h.s. for $u, v\in \plainH1(\Om)$
(below all integrals
are over $\Om$ unless indicated otherwise):
\begin{align*}
h_{N}^{(0)}(\om)& [\om^{-1}u, \om^{-1}v]
= \int \om^2\lu \BG \nabla(\om^{-1} u),
\overline{\nabla(\om^{-1} v)}\ru d\bx\\[0.2cm]
= &\ \int \om^2\lu
\BG(\om^{-1}\nabla u - \om^{-2}u\nabla\om ),
\om^{-1} \overline{\nabla v} -
\om^{-2}\overline{v}\nabla\om
\ru d\bx\\[0.2cm]
= &\ \int\lu \BG\nabla u, \overline{\nabla v}\ru d\bx
+ \int \om^{-2}\lu\BG\nabla\om,
\nabla\om\ru u
\overline{v}d\bx\\[0.2cm]
&\qquad - \int \om^{-1}
u\lu\BG\nabla\om, \overline{\nabla v}\ru d\bx
- \int \om^{-1} \lu\BG\nabla u,
\nabla\om\ru\overline{v} d\bx.
\end{align*}
Integrate the last integral by
parts,
remembering that $\BG\nabla\om\in\boldW1{q}, q > 1$,
and using the notation
$\mathbf n(\bx)$ for the exterior unit normal
to the boundary $\p\Om$ at the point $\bx$:
\begin{gather*}
- \int \om^{-1} \lu\BG\nabla u, \nabla\om\ru\overline{v} d\bx
= -\int_{\p\Om} \om^{-1}\lu \BG\nabla\om, \mathbf n\ru
u\overline{v}dS\\[0.2cm]
+ \int \om^{-1} u \lu\BG\nabla\om,\overline{\nabla v}\ru
d\bx
+ \int \om^{-1} u \overline{v}\lu\nabla,
\BG\nabla\ru \om  d\bx\\[0.2cm]
- \int \om^{-2}
\lu\BG\nabla\om, \nabla\om\ru u \overline{v} d\bx.
\end{gather*}
Substituting this in the
initial formula for the bilinear form,
we arrive at the relation
\begin{align*}
h_N^{(0)}(\om)[\om^{-1}u, \om^{-1} v]
= &\ h_N^{(0)}(1)[u, v]
+ \int\tilde V u \overline{v}d\bx
-\int_{\p\Om} \om^{-1}\lu \BG\nabla\om, \mathbf n\ru
u\overline{v}dS\\
= &\ h_N(1, \tilde V, \tilde{\bsig})[u, v].
\end{align*}
It remains to notice that
\begin{equation*}
\int_{\p\Om} \om^{-1}\lu \BG\nabla\om, \mathbf n\ru
u\overline{v}dS
= \sum_{j=1}^M \int_{\Sigma_j} \om^{-1}\lu \BG\nabla\om, \mathbf n\ru
u\overline{v}dS,
\end{equation*}
which completes the proof.
\end{proof}

\begin{cor}\label{om1:cor}
It suffices to prove Theorems \ref{rtwo:thm}
and \ref{dirneum1:thm} for $\om = 1$ and $\mu = 1$.
\end{cor}

\begin{proof}
To be definite, consider the
case of Theorem \ref{dirneum1:thm} only.
Suppose that it holds for $\om = \mu = 1$.
First we show that the operator
$H_{\aleph}(1, V, \mu),\ \aleph = D, N$, is absolutely
continuous if $\mu$ satisfies \eqref{weight:eq}.
According to Proposition
\ref{sc:prop} it suffices to check that
it has no point spectrum.
On the contrary, suppose that
$\l$ is an eigenvalue of $H_\aleph$
with an eigenfunction $u$.
Recall that by Lemma \ref{qforms:lem}
the quadratic form $h_\aleph$
is closed on the domain $\CD_\aleph$
\textsl{independent} of $\mu$.
Therefore the equality
\begin{equation*}
h_{\aleph}(1, V, \mu, \bsig)[u, v]
- \l\int_{\Om} u \overline{v} \mu d\bx = 0,
\forall v\in
\CD_\aleph,
\end{equation*}
implies that the point $\tilde\l = 0$ is
an eigenvalue of the operator
$\tilde H_{\aleph} = H_{\aleph}(1, V - \l\mu, 1, \bsig)$.
On the other hand,
the potential $V - \l\mu$ satisfies \eqref{ndcond:eq},
so that by Theorem \ref{dirneum1:thm}
with $\om=\mu=1$ the spectrum
of $\tilde H_{\aleph}$ is absolutely continuous.
This contradiction proves the claim.

Let us now remove the condition $\om = 1$.
Again, in view of Proposition
\ref{sc:prop} it suffices to show that
the operator $H_\aleph(\om, V, \mu, \bsig)$
has no point spectrum.
Assuming the contrary, we obtain from
\eqref{sandw:eq} that
for any eigenfunction $u$ of $H_\aleph(\om, V, \mu, \bsig)$
associated with an eigenvalue $\l$, the
function $\om u$ will be an eigenfunction of the operator
$\tilde H_\aleph
= H_\aleph(1, \tilde V - \l\om^{-2},
\mu, \tilde{\bsig})$
associated with the eigenvalue  $\tilde\l = 0$.
This contradicts the absolute continuity
of $\tilde H_\aleph$, which follows from Theorem
\ref{dirneum1:thm} with arbitrary $\mu$ and $\om = 1$,
in the same way as in the first part of the proof.
\end{proof}

Referring to this Corollary, from now on
\textbf{ we always  assume that $\om = \mu = 1$.}

\subsection{Floquet decomposition for $H_P$}
Let us describe the Floquet decomposition for the
operator $H_P$ for $\mu = 1, \om = 1$.
The other coefficients
are supposed to satisfy the minimal
smoothness conditions \eqref{ndcond:eq},
\eqref{s:eq}, \eqref{g:eq}.

Further we
follow almost word to word the papers
\cite{BS}, \cite{BirSol}.
To describe the spectral decomposition,
along with the operator
$H_P$ acting in $\plainL2(\CC) = \plainL2(1, \CC)$,
introduce the family of operators
$H_P(k), k\in\R$ defined by the quadratic form
\begin{equation}\label{qform1:eq}
h_P(k)[u] = \int_{\Ttwo} \biggl[ \lu\BG(\BD - \ba + k\be_1)u,
\overline{(\BD - \ba + k\be_1) u}\ru + V |u|^2 \biggr] d\bx +
\int_{\BXI} \bsig |u|^2 dS,
\end{equation}
on the domain $\tilde\CD_P = \plainH1(\Ttwo)$.
Similarly to \eqref{qforms:eq},
this form is bounded from below and closed.
Consequently, it uniquely defines
a self-adjoint operator, which we denote by
$H_P(k)$.
Since $H_P(k)$ has compact
resolvent for all $k$, its spectrum is purely  discrete.

To relate the family $H_P(k)$ to the operator $H_P$ defined
in the previous section, we identify the space $\plainL2(\CC)$
with the direct integral
\begin{equation*}
\GH = \int_0^1 \plainL2 (\Ttwo) dk
\end{equation*}
with the help of the unitary transform called
the
\textsl{Gelfand transform} $U:\plainL2(\CC)\to \GH$:
\begin{equation*}
(U f)(\bx, k) = e^{i 2\pi k x_1}
\sum_{m\in \Z} e^{-i 2\pi k m }
f(\bx + 2\pi m\be_1), \ \ \bx\in \CO.
\end{equation*}
The mapping $U$ can be easily shown to be unitary.
Then standard calculations (see e.g.
\cite{RS}, Ch. XIII, \S 16)
lead to the formula
\begin{equation}\label{direct:eq}
U H_P U^{-1} = \int_0^1 H_P(k) d k,
\end{equation}
in the sense that for any $f\in \CD_P$ one has
$(Uf)(\ \cdot\ , k)\in \tilde\CD_P$ a.a. $k\in \R$
and
\begin{equation*}
h_P[u, v] = \int_0^1 h_P(k)
[(Uu)(\ \cdot\ ,k), (Uv)(\ \cdot\ , k)] dk.
\end{equation*}
Extend the form $h_P(k)$ to all $k\in \mathbb C$.
An estimate similar to \eqref{infinite:eq}
shows that this form is closed and sectorial
for all $k\in \mathbb C$
(see \cite{Kato}, Theorem VI.3.9). Denote the
associated m-sectorial operator by $H_P(k)$.
Its resolvent is compact.
Since $H_P(k)^* = H_P(\overline{k})$,
the family $H_P(\ \cdot\  )$
forms a self-adjoint
analytic family of type (B) with compact resolvent
(see \cite{Kato}, Sect. VII.2).
By the Kato-Rellich theorem
(see \cite{Kato},
Theorem VII.3.9 and Remark
VII.4.22)
the eigenvalues $\l_j(k)$ of
$H_P(k)$ can be labeled  so as to make
them and associated eigenfunctions
analytic in $k\in\R$.

\subsection{Absolute continuity}
In the proof of Theorem \ref{dirneum1:thm}
we shall need the absolute continuity
of the operator
$H_P(1, \BB, \ba, V, 1, \bsig;\CC)$  with a
constant \textsl{diagonal}
matrix $\BB$:

\begin{thm}\label{per:thm}
Let $\om = \mu = 1$.
Suppose that the conditions
\eqref{ndcond:eq}, \eqref{s:eq}
are fulfilled and that  $\BG = \BB$ is a constant
\textbf{diagonal} matrix with positive entries.
Then the operator
$H_P(1, \BB, \ba,  V, 1, \bsig; \CC)$
is absolutely continuous.
\end{thm}

This result easily follows from the estimate
\eqref{low:eq} below,  established in \cite{BSS}.
Before stating this estimate properly, we
make the following remark.

In view of the analyticity properties
of $H_P(k)$, in order to prove the absolute continuity
of  the operator $H_P$,
it will suffice to show that the eigenvalues
$\l_j$ are not constant as
functions of $k\in\R$ (see \cite{RS}, Theorem XIII.86).
In its turn, this fact
will follow if we show
that for any $\l\in\R$ there is a value
$k\in\mathbb C$ such that the operator
$H_P(k) - \l$ is invertible
(see \cite{Kato}, Theorems VII.1.10, VII.1.9).
To ensure this property we
analyse the behaviour of the operator $H_P(k)$ for large
values of $\im k$.
The following Proposition is a special case of
Theorem 2.1 in \cite{BSS}.

\begin{prop} \label{lowerbound:thm}
Suppose that the conditions of Theorem
\ref{per:thm} are satisfied.
Then there exist
numbers $\b\in\R$
and $y_0>0$ such that
\begin{equation}\label{low:eq}
\|\bigl(H_P(\b + i y)\bigr)^{-1}\| \le C y^{-1},
\end{equation}
for all $y\ge y_0$.
\end{prop}

For any $\l\in\R$
Proposition \ref{lowerbound:thm}
guarantees the invertibility
of the operator
$H_P(\b + i y) - \l$ if $y$ is sufficiently
large.
This ensures that the eigenvalues
$\l_j(k)$ are not constant
in $k\in\R$ and therefore the
spectrum of the direct integral \eqref{direct:eq},
and hence
that of the operator
$H_P$, are absolutely continuous.
This completes the proof of Theorem \ref{per:thm}.

\section{Isothermal coordinates.
Proof of Theorem
\ref{rtwo:thm}}
\label{conf4:sect}

As it was already explained,
the proofs of Theorems \ref{rtwo:thm} and
\ref{dirneum1:thm} are based
on a reduction of the operator $H_\aleph(\BG)$ to
the canonical form, i.e. to the operator $H_\aleph(\BA)$
with a constant positive-definite
matrix $\BA$ (Recall again
that we may assume without loss of
generality that $\om = \mu = 1$).
This reduction is done using the so-called isothermal
coordinates.
The required properties of
this coordinate change are stated in Theorems
\ref{change:thm} and \ref{changeom:thm}. Their proof
is postponed until Section \ref{conf6:sect}.

We consider two mappings.
One is a homeomorphism of the entire plane
onto itself, and
the other is a homeomorphism of the periodic domain
$\Om$ onto the straight strip $\CS_{\dc}$.

We always denote $\BF = \sqrt{\BG}$ and
assume that
\begin{equation}\label{det1:eq}
\det\BG = 1,
\end{equation}
so that $\det\BF = 1$ as well.

\subsection{Change of variables}
Below we denote
\begin{equation*}
\BJ =
\begin{pmatrix}
0&-1\\
1&0
\end{pmatrix}.
\end{equation*}
The first theorem below describes
a suitable coordinate change in $\Rtwo$.

\begin{thm}\label{change:thm}
Let $\BG$ satisfy \eqref{g:eq} and \eqref{det1:eq}.
Then there exists a unique
homeomorphism
$\boldf = \boldf_{\Rtwo} = (f_1, f_2):\Rtwo\to\Rtwo$,\
$\boldf\in\boldH1_{\loc}(\R^2)$
such that
\begin{itemize}
\item[(i)]
$\boldf(0) = 0$,\ $ \boldf(2\pi\be_1) = 2\pi\be_1$,\
and $|\boldf(\bx)|\to \infty$ as $|\bx|\to\infty$;
\item[(ii)]
The components $f_1, f_2$ satisfy the equation
\begin{equation}\label{beltrami:eq}
\nabla f_2 = \BJ\BG \nabla f_1,\ \textup{a.a.}
\ \bx\in\R^2.
\end{equation}
\end{itemize}
Moreover, the map $\boldf$
possesses the following properties:
\begin{itemize}
\item[(iii)]
The Jacobian
$J_{\boldf}(\bx) = \det \bigl({\GD}\boldf(\bx)\bigr)$
is positive a.a. $\bx\in\R^2$.
The function $\boldf$ and its inverse $\boldf^{-1}$
both belong to $\boldW1\tau_{\loc}(\R^2)$ with some
$\tau >2$.
\item[(iv)]
For any $u\in \plainH1_{\loc}$ the usual chain rule holds:
\begin{equation*}
\nabla(u\circ\boldf) = [\GD\boldf]^T(\nabla u\circ\boldf),
\ \  \textup{a.a.}  \ \ \bx\in \R^2.
\end{equation*}
Moreover, for any $h\in \plainL1_{\loc}(\R^2)$
the function $h\circ\boldf$
belongs to
$\plainL1_{\loc}(J_\boldf, \R^2)$ and
\begin{equation*}
\int_{\boldf(\Om_1)} h(\by) d\by
= \int_{\Om_1} \bigl(h\circ\boldf\bigr)(\bx)
J_\boldf(\bx) d\bx
\end{equation*}
for any open bounded $\Om_1\subset\R^2$.
\item[(v)]
For $\bh_1 = 2\pi\be_1$ and some
linearly independent vector $\bh_2$
one has
\begin{equation*}
\boldf(\bx + 2\pi \bn) = \boldf(\bx) + n_1 \bh_1
+ n_2\bh_2,\ \forall \bx\in\R^2,
\end{equation*}
for all $\bn\in\Ztwo$.
\item[(vi)]
If
\eqref{holder:eq} is fulfilled, then
$\boldf\in \boldC{1+\a}(\Rtwo)$
and the Jacobian
satisfies the estimate
$J_{\boldf}(\bx)\ge c$ for all $\bx\in \Rtwo$;
\end{itemize}
\end{thm}

Note that the properties stated
in Theorem \ref{change:thm}(iv)
are certainly standard
for smooth, or even Lipshitz maps
$\boldf$. For homeomorphisms $\boldf$ of class
$\boldH1_{\loc}$ such ``standard'' results as
the chain rule or the change of variables
under the integral, are not obvious.
They require in addition that the functions
$\boldf$ and $\boldf^{-1}$ should map sets of measure
zero into sets of measure zero.
This property follows
from Theorem \ref{change:thm}(iii)
according to
\cite{GR}, Ch.5, \S 3.

Let us now state the appropriate result for an admissible domain
$\Om$ satisfying Condition \ref{boundary:cond}.
Recall that the boundary $\p\Om$ consists of two disjoint Jordan
curves $\ell_+$ and $\ell_-$. Without loss of generality we assume
that $\mathbf 0 \in \ell_-$. As defined in Sect. \ref{conf2:sect},
$Z\subset\p\Om$ is the set of the boundary
points where the smoothness of $\ell_+$, $\ell_-$ breaks down.

\begin{thm}\label{changeom:thm}
Suppose that an admissible domain
$\Om$ satisfies Condition \ref{boundary:cond}.
Let $\BG$ satisfy \eqref{g:eq}, \eqref{det1:eq}
and \eqref{holder:eq}.
Then there exists a unique
homeomorphism
$\boldf = \boldf_\Om = (f_1, f_2):\Om\to\CS = \CS_1$,
$\boldf\in\boldH1_{\loc}(\Om)$ such that
\begin{itemize}
\item[(i)]
$\boldf(\mathbf 0)
= \mathbf 0$,\ $ f_1(\bx) \to +\infty$ as $x_1\to+\infty$,
and  $f_1(\bx)\to -\infty$ as $x_1\to-\infty$;
\item[(ii)]
The components $f_1, f_2$ satisfy the equation
\eqref{beltrami:eq} for a.a. $\bx\in \Om$.
\end{itemize}
The map $\boldf$ satisfies the following properties:
\begin{itemize}
\item[(iii)]
For some number  $h > 0$ one has
\begin{equation*}
\boldf(\bx + 2\pi n \be_1) = \boldf(\bx) + 2\pi n h \be_1,\
\forall \bx\in\Om,
\end{equation*}
for all $n\in\Z$;
\item[(iv)]
$\boldf\in \boldC{1+\a}
\bigl(\overline{\Om}\setminus Z\bigr)$,
$\boldf^{-1}\in \boldC{1+\a}
\bigl(\overline{\CS}\setminus
\boldf(Z)\bigr)$,\ and the Jacobian
$J_{\boldf}(\bx)$ is positive everywhere in $\Om$.
Moreover,
for each $\bx_0\in Z$ there exist a number $\nu\in (0, 2]$ and
four non-degenerate H\"older-continuous
matrix-functions $\BM, \BT$
and $\boldsymbol\Phi$, $\boldsymbol\Psi$
with real-valued entries such that in the vicinity
of $\bx_0$ and $\bz_0 =
\boldf(\bx_0)\in \boldf(Z)$
one has the representations
\begin{equation}\label{singular:eq}
\begin{split}
\GD\boldf(\bx) = |\bx-\bx_0|^{\frac{1}{\nu}-1}
\boldsymbol\Phi\biggl(
\frac{\bx-\bx_0}{|\bx-\bx_0|}\biggr)
\BM(\bx),\\[0.3cm]
\GD \boldf^{-1}(\bz) = |\bz - \bz_0|^{\nu-1}
\boldsymbol\Psi\biggl(
\frac{\bz-\bz_0}{|\bz-\bz_0|}\biggr)
\BT(\bz).
\end{split}
\end{equation}
\end{itemize}
\end{thm}

The next Lemma establishes
some further properties of $\boldf$ that follow
from  Theorems  \ref{change:thm} and \ref{changeom:thm}.

\begin{lem}\label{change:lem}
Let the matrix $\BG$ and the map $\boldf$ be
as in Theorem \ref{change:thm} or
\ref{changeom:thm}.
Then
the following identities hold:
\begin{gather}
\nabla f_1 = - \BJ\BG \nabla f_2, \label{laplace2:eq}\\[0.2cm]
\langle\BF\nabla f_1, \BF\nabla f_2\rangle =
0,\label{orthog:eq}\\[0.2cm]
J_{\boldf} = |\BF\nabla f_1|^2 =
|\BF\nabla f_2|^2;\label{jacob:eq}\\
J_{\boldf}^{-1}\GD\boldf\  \BG\ \GD\boldf^T = \BI.
\label{changeg:eq}
\end{gather}
\end{lem}

\begin{proof}
As $\det\BF = \det\BG = 1$, a direct
calculation shows that
\begin{equation*}
\BJ = \BG\BJ\BG,\ \ \BJ = \BF\BJ \BF.
\end{equation*}
Noticing also that $\BJ^2 = -\BI$, from \eqref{beltrami:eq}
we obtain  \eqref{laplace2:eq} and the relation
\begin{equation*}
\BF\nabla f_2  = \BJ \BF\nabla f_1.
\end{equation*}
This implies the orthogonality \eqref{orthog:eq}
in view of the
obvious equality $\langle \BJ\bxi, \bxi\rangle = 0$,
$\forall \bxi\in \Rtwo$. It also yields the equality
$|\BF\nabla f_1|^2 = |\BF\nabla f_2|^2$.
The equality \eqref{changeg:eq} is a direct consequence
of \eqref{orthog:eq} and \eqref{jacob:eq}.

To prove \eqref{jacob:eq} compute the
Jacobian, using \eqref{laplace2:eq}:
\begin{align*}
J_{\boldf}(\bx) =
\dfrac{\p f_1}{\p x_1}\dfrac{\p f_2}{\p x_2}
- \dfrac{\p f_1}{\p x_2}
\dfrac{\p f_2}{\p x_1}
= &\ - \langle\nabla f_1, \BJ\nabla f_2\rangle\\
= & \ - \langle\nabla f_1, \BG \BJ\BG \nabla f_2\rangle
= \langle\nabla f_1, \BG \nabla f_1\rangle\\
= &\ |\BF\nabla f_1|^2 = |\BF\nabla f_2|^2.
\end{align*}
\end{proof}

\begin{rem}\label{qconform:rem}
By \eqref{jacob:eq} the norm $|\GD \boldf|$
can be estimated as follows:
$|\GD\boldf|^2\le K J_{\boldf}$
with some positive constant $K$ depending on the matrix $\BG$.
Recall that this inequality serves as a definition of the
so-called $K$-quasi-conformal maps (see \cite{GT}, Sect 12.1).

Similarly, $|(\GD\boldf)^{-1}|^2
= \bigl(|\GD\boldf| \BJ_{\boldf}^{-1}\bigr)^2
\le K J_{\boldf}^{-1}$.
\end{rem}

\subsection{Unitary transformation}
Notice that the mapping $\boldf = \boldf_{\Rtwo}$ constructed in
Theorem \ref{change:thm}
transforms the lattice $\BOLG = (2\pi\Z)^2$
into the lattice generated by the vectors $\bh_1, \bh_2$.
It is slightly more convenient to reduce this
lattice back to $\BOLG$ by applying the non-degenerate
linear transformation $\BR = \BR_{\Rtwo}: \Rtwo\to\Rtwo$
defined by the following two relations:
\begin{equation*}
\BR\bh_1 = 2\pi\be_1, \ \BR \bh_2 =  2\pi\be_2.
\end{equation*}
By Theorem \ref{change:thm}(v) the composite
mapping $\boldg = \BR\circ\boldf$ satisfies
\begin{equation}\label{gper:eq}
\boldg(\bx + 2\pi \bn) =
\boldg(\bx) +
2\pi n_1\be_1 + 2\pi n_2 \be_2,\ \
\forall \bx\in\R^2, \ \forall \bn\in\Z^2.
\end{equation}
Note that in view of
\eqref{changeg:eq}, we have
\begin{equation}\label{changeg1:eq}
\BA := J_{\boldg}^{-1}\GD\boldg\ \BG\  \GD\boldg^T
= (\det\BR)^{-1} \BR\BR^T.
\end{equation}
Similarly,
in Theorem \ref{changeom:thm} the one-dimensional
lattice $\bg_1$ is transformed
into the lattice $\bg_1^h = \{2\pi nh \be_1\}, n\in\Z$.
We rescale $\bg_1^h$ back to $\bg_1$  by applying the
transformation $\BR =\BR_{\Om}: \CS\to \CS$
defined as follows:
\begin{equation*}
\BR\be_2 = \be_2,\ \BR\be_1 = h^{-1}\be_1,
\end{equation*}
Then clearly, the mapping $\boldg = \BR\circ\boldf$ satisfies
the relation
\begin{equation}\label{gper1:eq}
\boldg(\bx + 2\pi n\be_1) = \boldg(\bx) + 2\pi n\be_1,\
\forall\bx\in \Om,\ \ \forall n\in\Z.
\end{equation}
Note also that the matrix $\BA$ in \eqref{changeg1:eq}
is diagonal in this case.
Depending on the context, below we use
either the notation $\boldg, \BA$ or $\boldg_\L,
\BA_\L$ where $\L$ assumes one
of the two values: $\Rtwo$ or $\Om$.

Denote
\begin{equation*}
\tilde\L =
\begin{cases}
\L, \ \ \textup{if}\ \ \L = \Rtwo
\ \ \textup{or}\ \ \CC;\\[0.2cm]
\CS, \ \ \textup{if}\ \ \L = \Om.
\end{cases}
\end{equation*}
In the case $\L = \Om$, under the conditions
of Theorem \ref{changeom:thm}
one can easily show that the set
$\tilde{\BSIG} = \boldg({\BSIG})$
 is again a system of curves.
More precisely,
if $\ell_j, j = 1, 2, \dots, N$ are the $\plainC1$-arcs
from Definition \ref{sigma:defn} (see also
end of Sect. \ref{conf2:sect}),
then each $\boldg(\ell_j)$ is
again a $\plainC1$-arc satisfying all
the required properties.
In the case of curves $\ell_j$ ending or starting
at the points of $Z$ this is done by a suitable
re-parametrisation, using \eqref{singular:eq}.
Recall that $\boldg({\BSIG})$ contains
the boundary of the strip $\CS$.

Introduce also the sets $\tilde\U$, $\tilde{\BXI}$
defined similarly to
\eqref{u:eq} and \eqref{xi:eq}  with $\tilde\L$ and
$\tilde{\BSIG}$ instead of $\L$ and $\BSIG$.
Note that according to Theorem \ref{changeom:thm}
the following is fulfilled
for the mappings $\boldg, \boldg^{-1}$
in the case $\L = \Om$:
\begin{equation}\label{w1t:eq}
\boldg\in\ \boldW1 \tau(\Om'),\
\boldg^{-1}\in \boldW1 \tau(\Om''),\ \
\left.\GD\boldg^{-1}\right|_{\tilde{\BXI}}\in
\boldL{\tau/2}(\tilde{\BXI})
\end{equation}
for any bounded domains $\Om'\subset \Om$, $\Om''\subset \CS$
and some $\tau > 2$.

Using  Theorem \ref{change:thm}(iii), (iv)
and Theorem \ref{changeom:thm} (iv),
it is easy to show that the operator
\begin{equation*}
(Su)(\bx) = u\bigl(\boldg^{-1}(\bx)\bigr),
\ u\in \plainL2(\Lambda),
\end{equation*}
is unitary from
$\plainL2(\L)$ onto the space
$\plainL2(\tilde\mu, \tilde\L)$ with the weight
\begin{equation*}
\tilde\mu(\bx) =
\biggl(J_{\boldg}\bigl(\boldg^{-1}(\bx)\bigr)\biggr)^{-1}
= J_{\boldg^{-1}}(\bx).
\end{equation*}
Denote
\begin{equation*}
\begin{cases}
\tilde \ba(\bx)
= \bigl(((\GD\boldg^T)^{-1}
\ba) \circ\boldg^{-1}\bigr)(\bx),\\[0.2cm]
\tilde
V(\bx) = \bigl( (J_{\boldg}^{-1}
V)\circ\boldg^{-1}\bigr)(\bx),
\end{cases}\ \ \bx\in \tilde\L,
\end{equation*}
and in the case $\L = \Om$ denote
\begin{gather*}
\tilde{{{{\bsig}}}}
= (\tilde\s_1, \tilde\s_2, \dots\tilde\s_N),\\[0.3cm]
\tilde\s_j(\bx) = |(\GD\boldg^{-1})(\bx) \mathbf t_j(\bx) |
(\s_j\circ\boldg^{-1})(\bx),\ \
\bx\in \tilde\Sigma_j,\ j = 1, 2, \dots, N,
\end{gather*}
where by $\mathbf t_j(\bx)$ we have
denoted the unit tangent vector
to $\tilde\Sigma_j$ at the point $\bx\in\tilde\Sigma_j$.

In the next Theorem, among other properties
of the map $S$ we show that
the coefficients
$\tilde V, \tilde\mu$, $\tilde\ba$ and
$\tilde{\bsig}$
satisfy the conditions of Theorems \ref{rtwo:thm}
or \ref{dirneum1:thm}.

\begin{thm}
Let $S$ be as defined above.
Then
\begin{itemize}
\item[(i)] Under the conditions of
Theorems \ref{rtwo:thm}
or \ref{dirneum1:thm} one has
$\tilde\mu\in\plainL {\tau/2}(\tilde\U)$,\
\begin{gather*}
\tilde V\in \plainL{\tilde p}(\tilde\U),\
\ \tilde p = \frac{p \tau}{2(p-1) + \tau},\\
\tilde \ba \in \boldL{\tilde s}(\tilde\U),\ \
\tilde s = \frac{s \tau}{s - 2 + \tau},\\
\tilde{\bsig}
\in \boldL{\tilde r}(\tilde{\BXI}),\ \
\tilde r = \frac{r \tau}{2 (r -1) + \tau}.
\end{gather*}
where $\tau > 2$ is as in Theorem
\ref{change:thm}(\textup{iii}) or \eqref{w1t:eq}.
The exponents $\tilde p$ and $\tilde s$ satisfy the inequalities
$\tilde p >1$, $\tilde s > 2$, $\tilde r > 1$;
\item[(ii)]
The map $S$ (resp. $S^{-1}$) is bounded as an
operator from $\plainH1(\L)$ to $\plainH1(\tilde \L)$
(resp. $\plainH1(\tilde\L)$ to $\plainH1(\L)$),
and from $\plainH1_0(\Om)$ to $\plainH1_0(\CS)$
(resp. $\plainH1_0(\CS)$ to $\plainH1_0(\Om)$).
\end{itemize}
\end{thm}

\begin{proof}
(i) The inequalities $\tilde p >1$, $\tilde s >2$
and $\tilde r >1$ immediately
result from the conditions $ p >1,\  s>2$,\
$r > 1,\  \tau>2$ by
inspection.

It follows from
Theorem \ref{change:thm}(iii) and (v) or
\eqref{w1t:eq},
that $\tilde\mu\in\plainL {\tau/2}(\tilde\U)$.
Let us prove that $\tilde V\in L^{\tilde p}(\tilde\U)$.
By H\"older's inequality for any bounded domain
$\Om_1\subset \R^2$ we have
\begin{multline*}
\int_{\Om_1} |J_\boldg^{-1} V\circ\boldg^{-1}(\bx)|^{\tilde p}
d\bx\\
\le \biggl[\int_{\Om_1}
|V\circ\boldg^{-1}(\bx)|^{\tilde p\b^{-1}} J_{\boldg^{-1}}(\bx)
d\bx\biggr]^{\b}
\biggl[\int_{\Om_1} \biggl(J_{\boldg^{-1}}(\bx)\biggr)^{(\tilde p -\b)
(1-\b)^{-1}}d\bx
\biggr]^{1-\b}
\end{multline*}
with $\b = \tilde p p^{-1}$.
Noticing that $(\tilde p -\b)
(1-\b)^{-1} = \tau/2$, and using Theorem \ref{change:thm}(iii),
(iv) or \eqref{w1t:eq},
we conclude that the r.h.s. of the last inequality
does not exceed
\begin{equation*}
C \biggl[\int_{\boldg^{-1}(\Om_1)} |V(\bx)|^p d\bx\biggr]^{\b}
\biggl[\int_{\Om_1}
 \biggl(J_{\boldg^{-1}}(\bx)\biggr)^{\tau/2}d\bx\biggr]^{1-\b}
<\infty.
\end{equation*}
To prove that $\ba\in\boldL{\tilde s}(\tilde\U)$
note that in view of Remark \ref{qconform:rem},
\begin{equation*}
\int_{\Om_1}|\tilde\ba(\bx)|^{\tilde s} d\bx
\le C \int_{\Om_1}
\biggl[ J_{\boldg^{-1}}(\bx)
\bigl|\bigl(\ba\circ\boldg^{-1}\bigr)
(\bx)\bigr|^2\biggr]^{\tilde s/2}
d\bx.
\end{equation*}
Now, repeating the argument used in the first part of the proof,
we arrive at the required property.

A similar calculation can be done for
the function $\bsig$.
Precisely, for any $\plainC1$-arc
$\ell\subset\tilde\Sigma_j$ we have
\begin{multline*}
\int_{\ell}
\biggl[
|(\GD\boldg^{-1})(\bx) \mathbf t_j(\bx) |\ \
|(\s_j\circ\boldg^{-1})(\bx)|\biggr]^{\tilde r}
d S\\
\le \biggl[\int_{\ell}
\biggl(
|(\s_j\circ\boldg^{-1})(\bx)|\biggr)^{\tilde r\g^{-1}}
|(\GD\boldg^{-1})(\bx) \mathbf t_j(\bx)|
dS\biggr]^{\g}\\
\biggl[\int_{\ell}
\biggl( |(\GD\boldg^{-1})(\bx)
\mathbf t_j(\bx)|\biggr)^{(\tilde r -\g)
(1-\g)^{-1}}d S
\biggr]^{1-\g}
\end{multline*}
with $\g = \tilde r r^{-1}$. Noticing that $(\tilde r -\g)
(1-\g)^{-1} = \tau/2$,
we conclude that the r.h.s. of the last inequality
does not exceed
\begin{equation*}
C \biggl[\int_{\boldg^{-1}(\ell)}
|\s_j(\bx)|^r d S\biggr]^{\g}
\biggl[\int_{\ell}
|\GD \boldg^{-1}(\bx)|^{\tau/2} d S\biggr]^{1-\g},
\end{equation*}
and therefore, by \eqref{w1t:eq} it is bounded.

(ii)
Let $u\in \plainH1(\L)$ and $v = Su$.
Let us first prove that
$\nabla v\in\plainL2(\tilde\L)$. Since the matrix
$\BA$ defined in \eqref{changeg1:eq}
is positive-definite,
changing the variables, we have
\begin{multline*}
c\|\nabla v\|^2
\le \int_{\tilde\L}
\lu\BA\nabla v, \overline{\nabla v}\ru d\bx
= \int_{\tilde\L} \lu \BG\nabla u,
\overline{\nabla u}\ru\circ\boldg^{-1} J_{\boldg^{-1}} d\bx\\
= \int_{\L} \lu\BG\nabla u,
\overline{\nabla u}\ru d\bx\le C\|\nabla u\|^2.
\end{multline*}
In the last inequality we have used
\eqref{g:eq}, and to secure the change
of variables in the case $\L = \R^2$ we refer to
Theorem \ref{change:thm}(iv).

Since
$\tilde\mu = J_{\boldg^{-1}}
\in \plainL {\tau/2}(\tilde\U)$
with $\tau>2$, in order to prove the
boundedness of the
operator $S:\plainH1(\L)\to \plainH1(\tilde\L)$
it remains to use the unitarity of
$S:\plainL2(\L)\to \plainL2(\tilde\mu, \tilde\L)$
and Lemma \ref{qforms:lem}.
Similarly for $S^{-1}$.

To show that $S$ maps $\plainH1_0(\Om)$ into
$\plainH1_0(\CS)$, it suffices to notice that
for any $u\in\plainC1_0(\Om)$
we have $u\circ\boldg^{-1}\in \plainC1_0(\CS)$.
\end{proof}

\begin{lem}\label{equivalence:lem}
Let the conditions \eqref{ndcond:eq}, \eqref{s:eq} be
fulfilled, let
$\boldg: \L\to\tilde\L$ be the homeomorphism
constructed above, and let $\BA$
be as defined in \eqref{changeg1:eq}.
Then the following unitary equivalence holds:
\begin{equation*}
S H_{\aleph}(\BG, \ba, V, 1,\bsig; \L) S^*
= H_{\aleph}(\BA,
\tilde\ba, \tilde V, \tilde\mu, \tilde{\bsig}; \tilde\L).
\end{equation*}
\end{lem}

\begin{proof}
Due to the properties of the mapping $S$ established above,
for $u, v\in \CD_\aleph(\tilde \L)$
we have $S^*u, S^*v\in \CD_\aleph(\L)$.
Therefore it suffices to show that
\begin{equation*}
h_\aleph(\BG, \ba, V, 1, \bsig; \L)[S^*u, S^*v]
= h_{\aleph}(\BA,
\tilde\ba, \tilde V, \tilde\mu, \tilde{\bsig}; \tilde\L)[u, v]
\end{equation*}
for any $u, v\in \CD_{\aleph}(\tilde\L)$.

Assume first that $V = 0$ and $\bsig = 0$.
Now denote  $\by = \boldg(\bx)$, and write the bilinear form
$h_{\aleph}(\BG, \ba, 0, 1, 0; \L)$:
\begin{align*}
h_{\aleph}& [S^*u, S^*v]\\
= &\ \int_{\L} \lu \BG(\bx)
(\nabla_{\bx} - i \ba(\bx))u\bigl(\boldg(\bx)\bigr),  \
\overline{(\nabla_{\bx} - i\ba(\bx))
v\bigl(\boldg(\bx)\bigr)}\ru d\bx\\
= &\ \int_{\L} J_{\boldg}(\bx)\\
&\ \times\lu \BA
\bigl(\nabla_{\by} - i(\BD\boldg^T)^{-1}\ba(\bx)\bigr)
u(\by),\
\overline{\bigl(\nabla_{\by} - i(\BD\boldg^T)^{-1}\ba(\bx)\bigr)
v(\by)}\ru d\bx.
\end{align*}
It remains to change the variables $\by = \boldg(\bx)$:
\begin{equation*}
h_\aleph[S^*u, S^*v]
= \int_{\tilde\L}
\lu \BA
\bigl(\nabla_{\by} - i\tilde\ba(\by)\bigr)
u(\by),\
\overline{\bigl(\nabla_{\by} - i\tilde\ba(\by)\bigr)
v(\by)}\ru d\by.
\end{equation*}
To take into account the terms with $V$ and $\bsig$ write
\begin{gather*}
\int_{\L} V(\bx) u\bigl(\boldg(\bx)\bigr)
\overline{v\bigl(\boldg(\bx)\bigr)} d\bx
+
\int_{\BSIG}\bsig(\bx) u\bigl(\boldg(\bx)\bigr)
\overline{v\bigl(\boldg(\bx)\bigr)} dS
\\
=
\int_{\tilde{\L}} \tilde V(\by) u(\by) \overline{v(\by)} d\by
+ \int_{\tilde{\BSIG}}
\tilde{\bsig}(\by) u(\by)
\overline{v(\by)} dS.
\end{gather*}
This calculation completes the proof.
\end{proof}

\begin{proof}[Proof of Theorem \ref{rtwo:thm}]
As it was explained in Remark \ref{constant:rem},
Theorem \ref{rtwo:thm} was established in \cite{BS2}
for constant matrices $\BG$ and $\om = \mu = 1$.
In view of the unitary equivalence established
in Lemma \ref{equivalence:lem}, and by Corollary
\ref{om1:cor}, Theorem \ref{rtwo:thm} for general $\BG$ and
$\om, \mu$ follows immediately from \cite{BS2}.
\end{proof}

As far as Theorem
\ref{dirneum1:thm} is concerned,
by Lemma \ref{equivalence:lem} and Corollary
\ref{om1:cor}, it will result from

\begin{thm}\label{dirneum:thm}
Suppose that $\Om = \CS_1$, $\om = \mu = 1$.
Let the conditions \eqref{ndcond:eq}, \eqref{s:eq}
be satisfied,
and let $\BG = \BB$ be a constant diagonal  matrix.
Then the spectra of the operators
$H_D(\BB, \ba, V, \bsig; \Om)$,
$H_N(\BB, \ba, V, \bsig; \Om)$
are absolutely continuous.
\end{thm}

\section{Proof of Theorems \ref{dirneum:thm} and
\ref{dirneum1:thm}}
\label{conf5:sect}

From now on we assume that the conditions of Theorem
\ref{dirneum:thm} are fulfilled.

\subsection{Extension operators}
The proof is based on a reduction of $H_D$
and  $H_N$ to a periodic
operator. To construct the appropriate
operator we reflect the
strip $\CS_1$ in the line $x_2 = 0$
and extend the functions $\ba,
V$ and $\s$ into the lower half of
the obtained domain. Precisely, denote
\begin{gather*}
\CS_u = \CS_1, \ \ \CS_l =
\{\bx\in \Rtwo: (x_1, -x_2)\in \CS_u\}, \\
\CS_0 = {\CS_u\cup\CS_l\cup\{\bx: x_2 = 0\}},
\end{gather*}
and define subspaces
\begin{equation*}
\plainL2_{\pm}(\CC) = \{u\in \plainL2(\CC):
u(x_1, x_2) = \pm u(x_1,
-x_2),\ \ \
\textup{a.a.}\  \bx\in \CS_0\}
\end{equation*}
of all even ($\plainL2_+$) and all
odd ($\plainL2_-$) functions from $\plainL2(\CC)$.
One easily  concludes by inspection that the
projections onto these subspaces are given by the formula
\begin{equation}\label{cp:eq}
\CP_\pm u = \dfrac{1}{2}
\bigl(u(x_1, x_2) \pm u(x_1, - x_2)\bigr).
\end{equation}
We shall also need \textsl{the extension operators}
$W_{\pm}:\plainL2(\CS_u)\to \plainL2_\pm(\CC)$.
For $\bx\in\CS_0$ they are
defined as follows:
\begin{equation*}
(W_\pm u)(x_1, x_2) =
\begin{cases}
u(x_1, x_2)/\sqrt{2},\ \bx\in \CS_u;\\[0.2cm]
\pm u(x_1, - x_2)/\sqrt{2},\ \bx\in \CS_l,
\end{cases}
\end{equation*}
and extended to $\Rtwo$ as $\bg_2$-periodic functions.
Using the formula
\begin{equation}\label{wadjoint:eq}
(W^{-1}_\pm u)(\bx) = \sqrt 2\  u(\bx), \bx\in\CS_u,\ \
\forall u\in \plainL2_{\pm}(\CC),
\end{equation}
one readily proves that
$W_{\pm}$ is unitary on $\plainL2(\CS_u)$.
To find out
how the functions
from the Sobolev spaces $\plainH1(\CS_u)$,
$\plainH1_0(\CS_u)$ transform under the
extensions $W_{\pm}$,
denote
\begin{equation*}
\plainH1_{\pm}(\CC) = \CP_{\pm}
\plainH1(\CC).
\end{equation*}
Given the explicit form \eqref{cp:eq}
of the operator $\CP_{\pm}$, it is easy to see that
\begin{equation}\label{intersect:eq}
\plainH1_{\pm}(\CC) = \plainH1(\CC) \cap
\plainL2_{\pm}(\CC).
\end{equation}
Now, observing that
\begin{equation}\label{zero:eq}
u(x_1, - \pi) = u(x_1, 0)
= u(x_1, \pi) = 0,\ \textup{a.a.}\
x_1\in \R,
\end{equation}
for all $u\in \plainH1_{-}(\CC)$, one sees that
\begin{equation}\label{extension:eq}
W_+\plainH1(\CS_u) = \plainH1_+(\CC),\ \ \
W_-\plainH1_0(\CS_u) = \plainH1_-(\CC).
\end{equation}

\subsection{Reduction}
Now we can describe the periodic operator associated with the
operators $H_D$ and $H_N$. Assume that the conditions of Theorem
\ref{dirneum:thm} are satisfied.
We begin with definition of the corresponding system of curves.
Define $\BSIG_u = (\Sigma_{1, u}, \Sigma_{2, u}, \dots, \Sigma_{N, u})$,
where $\Sigma_{j,u} = \Sigma_j, j = 1, 2, \dots, N$, and
\begin{gather*}
\BSIG_l = (\Sigma_{1, l}, \Sigma_{2, l}, \dots, \Sigma_{N, l}),\\
\Sigma_{j, l} =
\{\bx\in \overline{\CS_l}: (x_1, -x_2)\in \Sigma_{j, u}\},\\
\BSIG_0 = (\BSIG_u, \BSIG_l).
\end{gather*}
By Definition \ref{sigma:defn}  $\BSIG_0$
is a system of curves in $\CC$.
We emphasise that this system contains two copies of
the upper and lower boundaries of the strip $\CS_1$.

Define the $\bg_2$-periodic functions
$\bb, Q$ by applying the extension operators $W_{\pm}$ to
$\ba, V$ in the following way:
\begin{equation*}
Q = \sqrt2  \ W_{+}V,\ \
b_1 = \sqrt2 \ W_{+} a_1,\ \
b_2
= \sqrt2 \ W_{-} a_2.
\end{equation*}
It is clear that the functions
$Q, b_1$ are even, and $b_2$ is odd.
Clearly, the new coefficients
$\bb, Q$ satisfy \eqref{ndcond:eq} with $\U = \Ttwo$.
In a similar way we extend to $\BSIG_0$ the function $\bsig$,
i.e. let $\brho$ be the
function on $\BSIG_0$ defined as follows:
\begin{equation*}
\brho(\bx) =
\begin{cases}
\bsig(\bx),\ \ \bx\in \BSIG_u,\\
\bsig(x_1, -x_2),\ \bx\in \BSIG_l,
\end{cases}
\end{equation*}
and extended $\bg_2$-periodically.
Clearly, $\brho$ satisfies \eqref{s:eq} with
$\BXI = \BSIG_0/\bg_1$.

Now define the reference periodic operator
$H = H_P(1, \BB, \bb,  Q, 1, \brho; \CC)$.
Using the symmetry properties of
$\bb, Q, \brho$,
we decompose the operator
$H$ in the orthogonal sum associated
with the subspaces $\plainL2_\pm = \plainL2_\pm(\CC)$.

\begin{lem}\label{invariant:lem}
The subspaces $\plainL2_{\pm}$
are invariant subspaces of the operator $H_P$.
\end{lem}

\begin{proof}
We need to check that
\begin{equation}\label{orthog:eq1}
\CP_\pm D[h] = \plainL2_\pm \cap D[h],
\end{equation}
and that for any $u, v\in D[h]$
\begin{equation}\label{orthog:eq2}
h[u, v] = h[\CP_+u,\CP_+ v] + h[\CP_-u, \CP_-v].
\end{equation}
The equality \eqref{orthog:eq1}
follows from \eqref{intersect:eq}.
To prove \eqref{orthog:eq2} it will suffice to verify that
\begin{equation}\label{pm:eq}
h[\CP_+u,\CP_- v] = 0,\ \forall u,v\in \plainH 1(\CC).
\end{equation}
Write out the l.h.s. using the notation
$w_\pm = \CP_\pm w$:
\begin{equation*}
\sum_{l=1}^2 \bigl(b_{ll}
(D_l - b_l)u_+,
(D_l - b_l)v_- \bigr) + (Q u_+, v_-)
+ \int_{\BSIG_0} \brho u_+ \overline{v_-} dS.
\end{equation*}
Two last terms vanish since
$Q$ and $\brho$ are even functions.
In particular, the integrals over $\p\CS$
equal zero in view of \eqref{zero:eq}.
To handle the first term make
the following table using the
properties of the coefficients:
\begin{align*}
&(D_1 - b_1)u_+ \ \  \textup{even,}
&(D_1 - b_1) v_- \ \  & \textup{odd,}\\[0.2cm]
&(D_2 - b_2)u_+ \ \  \textup{odd,}
&(D_2 - b_2) v_- \ \   & \textup{even.}
\end{align*}Now it is easy to see that the first
term also vanishes, which implies \eqref{pm:eq}.
\end{proof}

Denote the parts of the operator $H$ in this orthogonal
decomposition by $H_+$ and $H_-$.
It can be shown that
$H_\pm$ are unique self-adjoint
operators associated with the
closed quadratic forms $h_{\pm}[\ \cdot\ ]$
obtained from $h[\ \cdot\ ]$ by
restricting
the domain $D[h]$ to $D[h_\pm] = \plainH 1_{\pm}(\CC)$.
The final step of the reduction of $H_D$ and $H_N$ to
$H$ is implemented in the following
Lemma:

\begin{lem}\label{unitary:lem}
Let the conditions of Theorem
\ref{dirneum:thm} be satisfied.
Then
\begin{equation}\label{grandfinal:eq}
W_- H_D W_-^* = H_-,\ \ W_+ H_N W_+^* = H_+.
\end{equation}
\end{lem}

\begin{proof}
By \eqref{extension:eq}
$W_+ \CD_N = D[h_+]$ and $W_-\CD_D = D[h_-]$,
and hence it suffices to show that the bilinear forms
of the operators in \eqref{grandfinal:eq},
coincide on the domains $D[h_-] = \plainH1_-(\CC)$ and
$D[h_+] = \plainH1_+(\CC)$ respectively.
Let $u_\pm, v_\pm\in \plainH1_\pm(\CC)$.
Using \eqref{wadjoint:eq} and referring again
to the symmetry properties of the coefficients,
one can write:
\begin{align*}
 & \bigl(\BB(\BD-\ba) W_{\pm}^* u_{\pm},
(\BD -\ba) W_{\pm}^* v_{\pm}\bigr)\\
&\qquad\qquad\qquad +
\bigl( V W_{\pm}^*u_{\pm}, W_{\pm}^*v_{\pm}\bigr)
+ \int_{\BSIG} \bsig W_{\pm}^*u_{\pm}
\overline{W_{\pm}^* v_{\pm}} dS\\
= &\ 2\int_{\CS_u}
\biggl[ \sum_{l=1}^2 b_{ll}(D_l-a_l)
u_{\pm} \overline{(D_l -a_l)v_\pm}
+ V u_{\pm}\overline{v_{\pm}}\biggr]d\bx
+ 2 \int_{\BSIG} \bsig u_{\pm}
\overline{v_{\pm}} dS\\
= &\ \int_{\CS_0}
\biggl[ \sum_{l=1}^2
b_{ll}(D_l - b_l)
u_\pm\overline{(D_l - b_l)v_\pm}
+ Q u_\pm\overline{v_\pm}\biggr]d\bx
+  \int_{\BSIG_0} \brho u_{\pm}
\overline{v_{\pm}} dS.
\end{align*}
This form coincides with the bilinear form of the
operator $H_\pm$.
\end{proof}

\begin{proof}[Proof of Theorem \ref{dirneum:thm}]
The coefficients
$\bb, Q$ satisfy the conditions
of Theorem \ref{per:thm}. Therefore,
the periodic operator $H$ is absolutely continuous, and
so are the orthogonal parts $H_+$ and $H_-$.
By virtue of Lemma \ref{unitary:lem}
the operators
$H_N$ and $H_D$ are unitarily equivalent to
$H_+, H_-$ and thus they
are also absolutely continuous,
as required.
\end{proof}

Theorem \ref{dirneum:thm}
combined with Corollary \ref{om1:cor}
and Lemma \ref{equivalence:lem} leads to Theorem \ref{dirneum1:thm}.

\section{Quasi-conformal maps}\label{conf6:sect}

In this section we prove
Theorems \ref{change:thm} and
\ref{changeom:thm}.
We are using the standard approach to second order
elliptic equations in dimension two, which consists
in passing to the complex variable and using the theory of
quasi-analytic functions
(see e.g. \cite{BJS}, \cite{Ve}).
Let us define
$z = x_1 + ix_2$, $f = f_1 + if_2$.

To begin with we notice that the equation
\eqref{beltrami:eq} for $f_1$, $f_2$
is equivalent to the \textsl{Beltrami} equation
(see \cite{Ve}):
\begin{equation}\label{complex:eq}
\p_{\bar z} f = q \p_z f,
\end{equation}
where
\begin{equation*}
\p_z = \frac{1}{2}\bigl(\p_{x_1} - i\p_{x_2}\bigr), \
\p_{\bar z} = \frac{1}{2}\bigl(\p_{x_1} + i\p_{x_2}\bigr)
\end{equation*}
with the complex-valued function
\begin{equation}\label{q:eq}
q =  \frac{-g_{12}+ i(1-g_{22})}
{g_{12} -i(g_{22} +1)}
\end{equation}
Note that $\|q\|_{\plainL\infty}$
is strictly less that $1$, and that the Jacobian
of $\boldf$ satisfies the relation:
\begin{equation*}
J_{\boldf} = |\p_z f|^2 - |\p_{\bar z} f|^2
= (1 - |q|^2)|\p_z f|^2.
\end{equation*}
We say that a continuous
function $f$ is $q$-quasi-conformal
(or quasi-conformal) if
$\p_zf\in\plainL2_{\loc}$
and $f$ satisfies the equation
\eqref{complex:eq}
for a.a. $z$.

\subsection{Mappings of $\Rtwo$}
The following theorem
is a one-to-one ``translation'' of Theorem
\ref{change:thm} into the language of
quasi-conformal mappings,
except for additional item (v).

\begin{thm}\label{periodicity:thm}
Let $q\in \plainL\infty(\mathbb C)$ be
a function such that $\|q\|_{\plainL\infty} < 1$.
Suppose that $q$ is periodic, that is,
\begin{equation*}
q(z) = q(z+2\pi) = q(z+2\pi i), \ \textup{a. a.}\
z\in \mathbb C.
\end{equation*}
Then there exists a unique $q$-quasi-conformal
homeomorphism $f$ of the complex plane onto itself
such that $f(0) = 0$, $f(2\pi) = 2\pi$,
and $f(\infty) = \infty$.

Moreover,
\begin{itemize}
\item[(i)]
$\p_z f\ne 0$ almost everywhere;
\item[(ii)]
There exists a number $\tau >2$ such that
$f, f^{-1}\in \plainW1{\tau}_{\loc}
(\mathbb C)$;
\item[(iii)] For any $u\in \plainH1_{\loc}(\mathbb C)$
the derivatives $\p_z (u\circ f)$ and $\p_{\bar z}(u\circ f)$
are found by the standard chain rule.
Also, for any $h\in \plainL1_{\loc}(\mathbb C)$
the function $h\circ\boldf$
belongs to
$\plainL1_{\loc}(J_\boldf, \mathbb C)$ and
\begin{equation*}
\int_{f(\Om_1)} h(\bx) d\bx
= \int_{\Om_1} \bigl(h\circ f\bigr)(\bx) J_\boldf(\bx)
d\bx
\end{equation*}
for any open bounded  $\Om_1\subset \mathbb C$;
\item[(iv)]
The mapping $f$ has the following periodicity property:
\begin{equation}\label{periodicity:eq}
f(z+2\pi n + 2 i\pi m) = f(z) + 2\pi n + \vark m,\
\forall z\in\mathbb C,\ \
\forall m,n,\in\Z,
\end{equation}
with some $\vark$ which has a non-zero imaginary part;

\item[(v)]
If $q(\overline{z}) = \overline{q(z)}$ almost everywhere,
then $\vark$ in \eqref{periodicity:eq} is purely imaginary:
$\vark = 2i\, \mbox{\rm Im} f(\pi i)$.
\item[(vi)]
If $q\in \plainC{\a}(\mathbb C),\ 0<\a<1$,
then $f\in \plainC{1+\a}(\mathbb C)$ and
$|\p_z f| > 0$ for all $z\in\mathbb C$.
\end{itemize}
\end{thm}

\begin{proof} By Proposition
\ref{global:prop} the mapping with the fixed three points $0,
2\pi, \infty$ exists, is unique, and satisfies properties
(i), (ii), (iii) and (vi). We need only to prove
\eqref{periodicity:eq} and (v).
Let $\z$ be either $2\pi$ or $2\pi i$.
Along with $f$ the function
\begin{equation*}
\tilde f(z) = 2\pi\frac{f(z+\z) - f(\z)}{f(2\pi+\z) - f(\z)}
\end{equation*}
is also a solution of the
equation \eqref{complex:eq},
due to the periodicity of
$q$.
Note that the denominator
here is not zero as the function
$f$ is a homeomorphism, and for the same reason
$f(\z) \not=0$ as well.
Notice also that
$\tilde f(0) = 0$, $\tilde f(2\pi) = 2\pi$ and
$\tilde f(\infty) = \infty$. By uniqueness
of such a solution,
$\tilde f = f$, or, which is the same,
\begin{equation*}
f(z+\z) = c_1 f(z) + c_2
\end{equation*}
with
\begin{equation*}
c_1 = \frac{f(2\pi+\z) - f(\z)}{2\pi},\ \
c_2 = f(\z).
\end{equation*}
Let us consider separately four possibilities:
$|c_1|<1$, $|c_1| > 1$, $c_1 = e^{i\t}, \t \in (0, 2\pi)$ and
$c_1 = 1$. We shall eliminate
the first three of them, thus proving
that $c_1 = 1$.

Case 1: $|c_1|<1$. For any integer $n\ge 1$ we have
\begin{equation*}
f(n\z) = c_2\sum_{k=0}^{n-1} c_1^k
= c_2\frac{1-c_1^n}{1-c_1},
\end{equation*}
and hence $f(n\z)\to c_2(1-c_1)^{-1}$ as $n\to\infty$.
This contradicts the fact that $f(\infty) = \infty$.

Case 2: $|c_1|>1$. The sought contradiction follows from
Case 1 by rewriting
\begin{equation*}
f(z-\z) = \frac{1}{c_1} f(z) -\frac{c_2}{c_1},
\end{equation*}
and noting that $|c_1^{-1}| <1$.

Case 3: $c_1 = e^{i\t},\ \t\in(0, 2\pi)$. Again, as in Case 1,
we have
\begin{equation*}
f(n\z) = c_2\frac{1-c_1^n}{1-c_1}.
\end{equation*}
The r.h.s. remains bounded as $n\to\infty$, which contradicts
the requirement that $f(\infty) = \infty$.

Consequently, the only possible option is $c_1 = 1$.

Note that $f(2\pi) = 2\pi$ by definition of $f$, so that
\eqref{periodicity:eq} with $m = 0$
is proved. Let us prove now that the imaginary part
of  $\vark = f(2\pi i)$ is non-zero. Suppose, on the contrary, that
$\im \vark = 0$.
Then for any integer $m$
one can find another integer $n = n(m)$ such that
$|2\pi n + \vark m|\le 2\pi$, so that the r.h.s.
of the equality
\begin{equation*}
f(2\pi n + 2\pi m i) = 2\pi n +  \vark m
\end{equation*}
remains bounded as $m\to\infty$. On the other hand
$|2\pi n + 2\pi m i|\to \infty$ as $m\to\infty$.
Again we get the same contradiction, which proves that
$\im \vark\not = 0$.

It is left to prove (v). Let $\varphi(z) :=
\overline{f(\overline{z})}$. Then $\varphi(0) = 0$,
$\varphi(2\pi) = 2\pi$ and $\varphi(\infty) = \infty$.
Further, it is easy to see
(see e.g. \cite{Ah}, Ch. I,
Sect. C, (1)) that
\begin{equation*}
\p_z (f(\overline{z})) = (\p_{\overline{z}} f)(\overline{z}) =
(q\p_z f)(\overline{z}) =
q(\overline{z}) (\p_z f)(\overline{z}) = q(\overline{z})
\p_{\overline{z}} (f(\overline{z})) .
\end{equation*}
Taking the complex conjugates of both sides and using
the equality
$q(\overline{z}) = \overline{q(z)}$, we obtain
$\p_{\overline{z}} \varphi (z) = q(z)\p_z \varphi(z)$.
By uniqueness we then conclude that $\varphi \equiv f$,
i.e. $\overline{f(\overline{z})} = f(z)$, $\forall z \in
\mathbb C$, and in particular, $f(-\pi i)
= \overline{f(\pi i)}$.
On the other hand \eqref{periodicity:eq} implies
\begin{equation*}
f(\pi i) = f(-\pi i + 2\pi i) = f(-\pi i) + \vark =
\overline{f(\pi i)} + \vark.
\end{equation*}
Hence $\vark = 2i\, \mbox{\rm Im} f(\pi i)$.
\end{proof}

Theorem \ref{change:thm} follows immediately.

We emphasise again that the crucial periodicity property
\eqref{periodicity:eq}
of the quasi-conformal map $f$ is
a direct consequence of the uniqueness in Proposition
\ref{global:prop}.
Besides,
we have included in Theorem \ref{periodicity:thm}
statement (v) which also follows from the uniqueness.
Using \eqref{q:eq} one easily sees that the condition
$q(\overline z) = \overline {q(z)}$ in Part (v) is
equivalent to the following symmetry
conditions on the matrix $\BG$:
\begin{gather*}
g_{jj}(x_1, x_2) = g_{jj}(x_1, -x_2), \ j = 1, 2,\\
g_{jl}(x_1, x_2)
= - g_{jl}(x_1, -x_2),\ j\not=l.
\end{gather*}
Then $\re\vark = 0$ means
that the isothermal change of variables transforms
the initial square lattice $(2\pi\Z)^2$
into another \textsl{orthogonal} lattice.
Although this observation is not needed in this paper, we
consider it to be worth mentioning.

\subsection{Mapping of the domain $\Om$}
The next theorem restates Theorem \ref{changeom:thm}
in the language of quasi-conformal mappings.

\begin{thm}\label{strip:thm}
Let $\Om$ be a domain as
in Theorem \ref{changeom:thm}, in particular
Condition \ref{boundary:cond} is fulfilled and
$0\in\ell_-$.
Let $q\in \plainC{\a}\bigl(\overline{\Om}\bigr)$
be a periodic function, i.e.
$q(z)  = q(z+2\pi),\ \forall z\in \Om$,
such that $\|q\|_{\plainL\infty}<1$.
Then there exists a unique $q$-quasi-conformal
homeomorphism $f$ of the domain
$\Om$ onto the strip $\CS_1$
such that $f(0) = 0$,
$f(-\infty) = -\infty$ and $f(+\infty) = +\infty$.
The map $f$ has the following properties:
\begin{itemize}
\item[(i)]
For all $z\in\Om$
\begin{equation}\label{periodstrip:eq}
f(z+2\pi n) = f(z) + \vark n,\ \forall n\in \Z,
\end{equation}
with some $\vark>0$;
\item[(iii)]
$f\in \plainC{1+\a}
\bigl(\overline{\Om}\setminus Z\bigr)$, \
$f^{-1}\in \plainC{1+\a}
\bigl(\overline{\CS}\setminus f(Z)\bigr)$
and $|\p_z f| > 0$ everywhere in $\Om$.
Moreover,
for each $\bz_0\in Z$ there exist a number $\nu\in (0, 2]$ and
four H\"older-continuous functions $M, T$ and $\Phi, \Psi$,
separated from zero,
such that in the vicinity
of $z_0$ and $\z_0 =
f(z_0)\in f(Z)$
one has the representations
\begin{gather*}
\p_z f(z) = |z - z_0|^{\frac{1}{\nu}-1}
\Phi\bigl(\arg(z - z_0)\bigr) M(z),\\
\p_\z f^{-1}(\z) = |\z - \z_0|^{\nu-1}
\Psi\bigl(\arg(\z - \z_0)\bigr)T(\z).
\end{gather*}
\end{itemize}
\end{thm}

\begin{proof}
Pick the following accessible boundary
points (see e.g. \cite{Gol}, Ch II, \S 3 for definition)
of $\p\Om$: $-\infty, +\infty, 0$.
Then by Theorem \ref{zhenya:thm}(i)
combined with Remark \ref{strip:rem},
there exists a uniquely defined
quasi-conformal homeomorphism $f:\Om \to \CS_1$ which
preserves these boundary points.
Moreover, the smoothness properties
required in (iii)
follow from
Theorem \ref{zhenya:thm}(ii) and (iii).
It remains to prove
\eqref{periodstrip:eq}.

In view of the periodicity of the domain $\Om$
the mapping $\tilde f(z) = f(z+2\pi)$ is also
a  quasi-conformal homeomorphism of $\Om$ on the straight
strip $\CS_1$.
Note that $\tilde f$ sends $-\infty$ and $\infty$
into themselves, and the point $0$
into $z_0 = f(2\pi)$.
Note that $z_0$ is real, since
it lies on the lower portion of
the boundary of $\CS_1$, i.e. on
the straight horizontal line $\im z = 0$.
Consequently, the composition function
\begin{equation*}
h(z) = \tilde f\bigl(f^{-1}(z)\bigr),
\end{equation*}
defined on $\CS_1$, is a conformal homeomorphism
of $\CS_1$ onto itself
(see \cite{BJS}, Part II, Sect. 6.2)
acting in such a way that
$-\infty$, $+\infty$ are preserved and
$h(0) = z_0$.
It is easy to see that the function
$\tilde h(z) = z + z_0$
satisfies the same conditions.
On the other hand,
such a conformal mapping is
unique (see e.g. \cite{Gol}, Ch. II, \S 3, Th. 6).
Therefore
$h(z) = \tilde h(z)$ and hence
\begin{equation*}
\tilde f(z) = f(z) + \vark,\ \vark
= \overline{\vark} = z_0.
\end{equation*}
It is clear that $\vark\not=0$, for otherwise
the function $f$ would have remained bounded
as $|z|\to\infty$, which would contradict
the assumption that $f$ sends $\pm\infty$
into itself. For the same reason
$\vark > 0$, since otherwise $+\infty$ and $-\infty$
would exchange their places under the mapping $f$.
\end{proof}

For conformal mappings the above argument can be found in
\cite{Gol}, Ch. V, \S 1.

Now Theorem
\ref{changeom:thm} follows from Theorem \ref{strip:thm}.

\subsection{Bilipschitz mappings of the domain $\Om$}
In this subsection we address a question
which has no direct effect on the results
of the paper, but is nevertheless
natural and important.
If the domain $\Om$
has no corners or peaks, then according to
Theorem \ref{strip:thm} the homeomorphism
$f:\Om\to \CS_1$ is $\plainC{1+\a}$-smooth. This
is guaranteed by the initially assumed
$\plainC{1+\a}$-smoothness
of the boundary $\p\Om$ (see Condition \ref{boundary:cond}) and
$\plainC{\a}$-smoothness of the
matrix $\BG$ (see \eqref{holder:eq}).
In the presence of corners or/and peaks one or both
of the derivatives
$\p_z f$, $\p_z f^{-1}$ are unbounded.
In particular, the homeomorphism $f$ may not
be Lipshitz even if the boundary $\p\Om$ is Lipshitz.
Therefore it is legitimate
to ask whether a domain $\Om$ with a
Lipshitz boundary $\p\Om$
admits a periodic bilipshitz map $\Om\to \CS_1$.

To state the question in a precise form recall that
a mapping $F$ from a metric space $(X_1, d_1)$ into
a metric space $(X_2, d_2)$ is called {\it bilipschitz} if
there exists a constant $M > 0$ such that
\begin{equation*}
d_1(x, y)/M \le d_2(F(x), F(y)) \le M d_1(x, y) , \ \ \forall
x, y \in X_1 .
\end{equation*}
We say that a curve $\ell\in\mathbb C$ is Lipshitz if it is
a bilipshitz image of $\R$.
Note by passing that it is easy to give an
intrinsic characterisation
of a Lipschitz curve. First of all,
it is clear that a  Lipschitz curve
is Jordan and locally rectifiable.
Conversely, let $\ell : \R \to \mathbb C$ be a
Jordan locally rectifiable curve such that
$|\ell(t)| \to \infty$ as $t \to \pm\infty$.
Using  the arclength parametrisation, one can easily
show that $\ell$ is a Lipschitz curve if and only if
it is an chord-arc curve, i.e. if there exists a constant
$K \ge 1$ such that the length of the subarc of $\ell$
joining any two points is bounded by $K$ times the
distance between them.

Assume that the boundary of a periodic domain
$\Om$ consists of two
disjoint Lipshitz curves $\ell_+, \ell_-$.
Our objective is to find a
bilipschitz mapping $F$ of $\CS_1$ onto
$\Om$ such that
\begin{equation}\label{periodlip:eq}
F(z + 2\pi) = F(z) + 2\pi , \ \  \forall z \in \CS_1 .
\end{equation}
It is evident that the bilipshitz regularity of the
curves $\ell_{\pm}$ is necessary for the existence
of such a mapping.
Let us convince ourselves that this condition is also
sufficient:

\begin{thm}\label{lipstrip:thm}
Let $\Om$ be a simply connected periodic domain with
a boundary consisting of two disjoint Lipschitz
curves $\ell_+$ and $\ell_-$. Then there exists
a bilipschitz mapping $F$  of  $\CS_1$ onto $\Om$
satisfying \eqref{periodlip:eq}.
\end{thm}

\begin{proof}
Let $\varphi: \CS_1 \to \Om$
be a conformal homeomorphism,
mapping $\pm\infty$ onto $\pm\infty$,
$0$ onto a given
point of $\ell_{-}$ and
such that $\varphi(\z + \vark) = \varphi(\z) + 2\pi$,
$\forall\z\in \CS_1$
with some $\vark > 0$.
The existence of such a map follows from \cite{Gol}, Ch. V, \S 1.
Since $\ell_{\pm}$ are
Jordan locally rectifiable curves,
$\varphi$ can be extended
to a homeomorphism of the closure of $\CS_1$ onto
the closure of $\Om$
(see \cite{Gol}, Ch. II, \S 3,
Theorem 4),
for almost all $x_0 \in \R$ the finite limits
\begin{equation*}
\lim_{\tau \to 0+} \varphi'(x_0 +i\tau) \not = 0, \ \
\lim_{\tau \to 0+} \varphi'(x_0 +i\pi-i\tau) \not = 0
\end{equation*}
exist,
and $\varphi$ is conformal at $x_0$ and $x_0 + i\pi$
(see \cite{Gol}, Ch. X, \S 1, Theorems 1 and 3).
The latter implies that the curves
$\ell_{\pm}$ meet the arcs
$\ell_0 = \varphi([x_0, x_0+i\pi])$ and
$\ell_1 = \varphi([x_0 + \vark, x_0 + \vark +i\pi])$
at the right angle.
Now let us define a closed ``fundamental domain'' of
the strip $\CS_1$:
\begin{equation*}
\Delta := \{x + iy : \ x_0 \le x \le x_0 + \vark, \
0 \le y \le \pi\}.
\end{equation*}
The boundary of $\Delta$ consists
of four segments
\begin{gather*}
\g_{-} = [x_0,  x_0 + \vark],\
\g_{+} = [x_0 + i\pi, x_0 + \vark + i\pi],\\
\g_0 = [x_0, x_0 + i\pi],\ \g_1 = [x_0, x_0 + \vark + i\pi].
\end{gather*}
The next step is to define a bilipshitz map $F_0$
from $\p\Delta$ onto $\varphi(\p\Delta) = \p\varphi(\Delta)$.
Let $\tilde{\ell}_{\pm} := \varphi(\g_{\pm})\subset\ell_{\pm}$.
Since $\ell_{\pm}$ are Lipshitz curves,
there exists a bilipshitz homeomorphism
$F_0: \g_{\pm}\to \tilde{\ell}_{\pm}$.
Further, set  $F_0(z) = \varphi(z)$ if $z\in \g_0$ or $\g_1$.
In view of the periodicity of $\varphi$, we have
\begin{equation}\label{per:eq}
F_0(z + \vark)  = F_0(z) + 2\pi,\ z\in\g_0.
\end{equation}
Since the arcs $\tilde{\ell}_{+}$, $\tilde{\ell}_{-}$ and
$\ell_0$, $\ell_1$
meet at the right angle,
it is easy to see that $F_0$ is a
bilipschitz mapping
 of $\p\Delta$
onto $\p\varphi(\Delta) =
\varphi(\p\Delta)$.
Then it follows from \cite{Tu} (see also \cite{Mac}) that
$F_0$ can be extended to a bilipschitz homeomorphism
of $\mathbb C$ onto itself. This defines a
bilipschitz homeomorphism
$F_0 : \Delta \to \varphi(\Delta)$.
By virtue of \eqref{per:eq}
it is now straightforward to see
that the extension $F_1$ of $F_0$, defined by
\begin{equation*}
F_1(\zeta + n\vark) := F_0(\zeta) + 2\pi n , \ \ \zeta \in
\Delta, \  n \in \Z,
\end{equation*}
is a bilipschitz mapping of $\CS_1$ onto $\Om$ such that
$F_1(\z + \vark) = F_1(\z) + 2\pi$,
$\forall \z \in
\CS_1$.
It remains to define $F$ by
\begin{equation*}
F(x + iy) := F_1\left(\frac\vark{2\pi}\, x + iy\right).
\end{equation*}
\end{proof}

\section{General properties of quasi-conformal mappings}
\label{conf7:sect}

We begin with describing a quasi-conformal homeomorphism
of the complex plane onto itself.
The following general result can be
found in \cite{Ah}, \cite{AB}, \cite{BJS}
Part II, Ch. 6, or \cite{Ve}, Ch. 2.

\begin{prop}\label{global:prop}
Let $q\in \plainL\infty(\Rtwo)$ be a function such that
$\|q\|_{\plainL\infty}< 1$.
Then there exists a unique
quasi-conformal homeomorphism $f$ of the complex plane such that
$f(0) = 0$, $f(2\pi) = 2\pi$ and $f(\infty) = \infty$.
Moreover,
\begin{itemize}
\item[(i)]
The conclusions (i), (ii), (iii) of Theorem
\ref{periodicity:thm} hold;
\item[(ii)]
If $q\in \plainC{\a},\ 0<\a<1$,
in a neighbourhood of some point $z_0$, then
$f\in \plainC{1+\a}$ and $|\p_z f| > 0$
in this neighbourhood.
\end{itemize}

\end{prop}

The next theorem is a ``quasi-conformal''
version of Riemann's
mapping theorem and
the results on the boundary behaviour of conformal
mappings.
The results
collected in this Theorem probably
are not new, but we have not been able to find them
stated in the form convenient for us.

\begin{thm}\label{zhenya:thm}
Let $\Om \subset \mathbb C$ be a simply
connected open domain with more than one boundary point,
such that all
points of $\p\Om$ are
accessible. Suppose
$q\in \plainL{\infty}(\Om)$ and
$\|q\|_{\plainL\infty} < 1$. Then
\begin{itemize}
\item[(i)] There exists a unique
$q$--quasiconformal map $w$ of $\Om$ onto the unit
disk $\CD := \{\zeta \in \mathbb C : \ |\zeta| < 1\}$
which maps three given
points
$z_1, z_2, z_3 \in \p\Om$ indexed in
order of their occurrence as one proceeds
in the positive direction along
$\p\Om$ (see \cite{Gol}, Ch. II, \S 3)
onto three given points
$\zeta_1, \zeta_2, \zeta_3 \in \p\CD$
similarly indexed.
This map and its inverse $f^{-1}$ belong to
$\plainW1\tau_{\loc}(\Om)$
for some $\tau > 2$.
The map $f$ defines a homeomorphism of the
compactification $\widehat{\Om}$ of $\Om$
by the accessible points
of $\p\Om$ onto the closed unit disk.
\item[(ii)]
Let $m \in \mathbb N\cup\{ 0\}$, $0 < \a < 1$,
$z_0 \in \widehat{\Om}$,
and let $q$ be
$\plainC{m + \a}$--smooth in a
neighbourhood of $z_0$. Suppose also
that if $z_0 \not\in \Om$
then the intersection of $\p\Om$
with a neighbourhood of $z_0$ is a
$\plainC{m + 1 + \a}$--smooth
Jordan curve. Then $f$ is
$\plainC{m + 1 + \a}$--smooth in a
neighbourhood of $z_0$ and
$|\p_z f(z_0)|>0$.
\item[(iii)] Let $z_0 \in \p\Om$ and suppose that
the intersection of
$\p\Om$ with a neighbourhood of $z_0$ is a piece-wise
$\plainC{1 + \a}$--smooth Jordan curve,
$0 < \a< 1$, with the only
angular point at $z_0$, which is not an outward cusp,
i.e. the interior with respect to
$\Om$ angle at $z_0$ is nonzero.
Let $\g : [t_1, t_2] \to \mathbb C$, $\g(t_0) = z_0$
be a parametrization of this
curve in the positive direction. Suppose also
that $q$ is $\plainC\a$--H\"older
continuous in a neighbourhood of $z_0$.
Then $f$ and its inverse $g := f^{-1}$
satisfy the following conditions
\begin{equation}\label{(2)}
\p_z f(z) = \left((z - z_0) + q(z_0)(\bar{z} -
\bar{z}_0)\right)^{\frac1{\nu} - 1}F(z) , \ \ \
\end{equation}
\begin{equation}\label{(3)}
\p_\z g(\z) = (\z - \z_0)^{\nu - 1}G(\z) , \ \ \
\p_{\bar{\z}} g(\z) = -q(g(\z))
\overline{\p_\z g(\z)} ,
\end{equation}
where $\z_0 = f(z_0)$,
\begin{equation}\label{(4)}
\nu = \frac1\pi \arg\left\{-\frac{\g'(t_0 - 0) +
q(z_0)\overline{\g'(t_0 - 0)}}{\g'(t_0 + 0) +
q(z_0)\overline{\g'(t_0 + 0)}}\right\} \in (0, 2],
\end{equation}
and $F$, $G$ are H\"older continuous and nowhere zero in
some neighbourhoods of $z_0$ and $\z_0$ respectively.
\end{itemize}
\end{thm}

\begin{rem}\label{strip:rem}
The theorem above
uses the disk $\CD$ as a ``target'' domain.
This choice was made only for definiteness and convenience
of the proof.
One can easily restate Theorem \ref{zhenya:thm}
choosing other simply connected domains
as targets.
In particular, making obvious modifications,
$\CD$ can be
replaced with the strip $\CS_1$.
This can be done by mapping $\CD$ onto $\CS_1$ using a
standard conformal map, and
noticing that composition of a conformal
and a $q$-quasi-conformal
mapping is again $q$-quasi-conformal
(see \cite{BJS}, Part II, \S 6.2)
\end{rem}

\begin{proof}[Proof of Theorem \ref{zhenya:thm}]
\underline{Step I} Take an arbitrary
extension $q_0 \in \plainL\infty(\mathbb C)$ of $q$
such that $\|q_0\|_{\plainL\infty} < 1$.
There exists a $q_0$--quasiconformal
homeomorphism $w : \mathbb C \to \mathbb C$
which belongs to
$\plainW1{\tau}_{\loc}(\mathbb C)$
with a $\tau >2$, together with its
inverse  (see
Proposition \ref{global:prop}).
Let $\Om_0 := w(\Om)$, $z_k' := w(z_k)$,
$k = 1, 2, 3$. It is clear
that $\Om_0$ is a simply connected
domain and all points of
$\p\Om_0$ are accessible.
Hence there exists a unique conformal
map $\psi$ of $\Om_0$ onto $\CD$ which maps
$z'_1, z'_2, z'_3$ onto
$\z_1, \z_2, \z_3$ (see, e.g., \cite{Gol}, Ch. II,
\S 3, Theorem 6). It is not difficult to see that
$f := \psi\circ w : \Om \to \CD$ is a $q$--quasiconformal map
(see \cite{BJS}, Part II, \S 6.2)
having all the properties announced in (i).

\underline{Step II} Let $f_1 : \Om \to \CD$
be an arbitrary  $q$-quasiconformal
homeomorphism mapping $z_1, z_2, z_3$
onto  $\z_1, \z_2, \z_3$. Then
$f_1\circ f^{-1} :\CD \to \CD$ is
an analytic (see \cite{BJS}, Part II, \S 6.2)
homeomorphism and  hence a
conformal automorphism
(see \cite{Gol}, Ch. II, \S 1)
with  fixed points $\z_1, \z_2, \z_3$.
By the uniqueness result  for conformal
maps we have $f_1\circ f^{-1}(z)
\equiv z$, i.e.  $f_1 = f$.
This proves uniqueness and shows that the map $f$
constructed above
{\bf does not depend on the choice of an
extension $\mathbf{q_0}$}.

\underline{Step III}
Under the conditions of (ii)
there exists
an extension $q_0$ which is $\plainC{m + \a}$-smooth
in a neighbourhood of $z_0$.
Then it follows from
\cite{Ve},
Theorem 2.9 and the proof of Theorem 2.12,
that $w$ from Step I is
$\plainC{m + 1 + \a}$-smooth
in a neighbourhood of $z_0$ and
\begin{equation}\label{(5)}
|\p_z w(z_0)|^2 - |\p_{\bar z} w(z_0)|^2 > 0 .
\end{equation}
So, (ii) will follow if we
prove that the conformal map $\psi$
from Step I is $\plainC{m + 1 + \a}$-smooth
in a neighbourhood of
$w(z_0)$ and $\psi'(w(z_0)) \not= 0$.
We need to do this only if
$w(z_0) \in \p\Om_0$, i.e. $z_0 \in  \p\Om$.

\underline{Step IV}
It follows from Step III that under the
conditions of Part (ii), the intersection of
$\p\Om_0$ with a neighbourhood of $w(z_0)$ is a
$\plainC{m + 1 + \a}$-smooth Jordan curve.
Let us take a simply connected subdomain
$\Om_1 \subset \Om_0$ with a
$\plainC{m + 1 + \a}$-smooth boundary
such that $\G := \p\Om_1\bigcap\p\Om_0$ is a
$\plainC{m + 1 + \a}$-smooth subarc of
the above Jordan curve,
$w(z_0)$ belongs to $\G$ and is different
from its endpoints.
By the uniqueness result for
conformal maps we can choose
conformal maps $\psi_1 : \Om_1 \to \CD$
and $\varphi : \CD \to \psi (\Om_1)$ so that
$$
\psi = \varphi\circ\psi_1 \ \ \mbox{in} \ \ \Om_1 .
$$
The map $\psi_1$ is $\plainC{m + 1 + \a}$-smooth
in the closure of $\Om_1$ and
$\psi'_1(w(z_0)) \not= 0$
(see \cite{Pom}, Theorems 3.5, 3.6)

The mapping $\varphi$ maps an open circular arc containing
$\psi_1(w(z_0))$ onto a
circular arc.
Hence $\varphi$ is holomorphic  in a neighbourhood of
$\psi_1(w(z_0))$ and $\varphi'(\psi_1(w(z_0))) \not= 0$
(see \cite{Gol}, Ch. II, \S 3, Theorem 5).
Therefore $\psi$ is $\plainC{m + 1 + \a}$-smooth
in a neighbourhood of  $w(z_0)$ and
$\psi'(w(z_0)) \not= 0$. Note by passing that this is in fact
a variant of  Kellogg's theorem
(cf., e.g. \cite{LS}, Sect. 29 or \cite{La}, Ch. II,
\S 1, Theorem 1).

The proof of (ii) is now completed.

\underline{Step V}
The proof of (iii) is similar to that of (ii).
Let $q_0$ be an extension of $q$ which is
$\plainC\a$-H\"older
continuous in a neighbourhood of $z_0$.
Then similarly to Step III, $w$ from Step I
is $\plainC{1 + \a}$-smooth in a
neighbourhood of $z_0$ and it satisfies
\eqref{(5)}. Therefore the intersection of
$\p\Om_0$ with a neighbourhood of $w(z_0)$
is a piece-wise
$\plainC{1 + \a}$-smooth Jordan
curve with the only singular
point at $w(z_0)$.
This curve is parameterized by
$[t_1, t_2] \ni t \mapsto w(\g(t)) \in \mathbb C$.
The equality $\p_{\bar z} w(z)
= q(z)\p_z w(z)$ implies
\begin{gather*}
\frac{d w(\g(t_0 \pm 0))}{dt} =
\p_z w(z_0)\g'(t_0 \pm 0) +
\p_{\bar z} w(z_0) \overline{\g'(t_0 \pm 0)}
\\
= (\g'(t_0 \pm 0) + q(z_0)\overline{\g'(t_0 \pm 0)})
\p_z w(z_0).
\end{gather*}
Therefore the interior with respect
to $\Om_0 = w(\Om)$ angle at
$w(z_0)$ equals $\pi\nu$, where
$\nu$ is given by \eqref{(4)}.

\underline{Step VI}
Now we need to investigate the properties of the
conformal map $\psi : \Om_0 \to \CD$.
Using the argument from Step IV
we can reduce this to the study of a conformal map
$\psi_1 : \Om_1 \to \CD$, where $\Om_1$ is a
simply connected domain with a piece-wise
$\plainC{1 + \a}$--smooth boundary
and the only singular point
$w(z_0)$, where the interior with respect
to $\Om_1$ angle equals $\pi\nu>0$.
Applying Warschawski's theorem (see e.g.
\cite{Ve}, Theorem 1.9,
and \cite{KKP}, Ch. 3, \S 3)
we obtain that $\psi$ and
its inverse $\eta := \psi^{-1}$ satisfy the conditions
$$
\psi'(w) = (w - w(z_0))^{\frac1{\nu} - 1}\Psi(w) , \ \ \
\eta'(\z) = (\z - \z_0)^{\nu - 1}E(\z) ,
\
\z_0 = \psi(w(z_0)),
$$
where $\Psi$ and $E$ are H\"older
continuous and nowhere zero in
some neighbourhoods of $w(z_0)$ and
$\z_0$ correspondingly.
Now (iii) follows from the formulae $f = \psi\circ w$,
$g = w^{-1}\circ\eta$ and
\begin{gather*}
w(z) - w(z_0) = \left(\p_z w(z_0)(z - z_0) +
\p_{\bar z} w(z_0)(\bar{z} - \bar{z}_0)\right)
\om(z, z_0)\\
=\left((z - z_0) + q(z_0)(\bar{z} -
\bar{z}_0)\right)\p_z w(z_0)\,\om(z, z_0),
\end{gather*}
where $\om(\cdot, z_0)$ is H\"older continuous in a
neighbourhood of $z_0$ and $\om(z_0, z_0) = 1$. The second
equality in \eqref{(3)} follows
from \cite{BJS}, Part II,  Ch. 6, Appendix, Theorem 3(iv)
(see also \cite{Ah}, Ch. I, Sect. C).
\end{proof}

\section{Acknowledgements}

The authors are thankful to
M. Birman, T. Suslina, and R. Shterenberg
for fruitful discussions.
R. Shterenberg's comments allowed us
to remove unnecessary restrictions on the
coefficients and thereby simplify the statements
and proofs of the main results.
We are also grateful to T. Suslina and R. Shterenberg
for making available to us the results
to be published in paper
\cite{ShSus}.

\bibliographystyle{amsplain}

\providecommand{\bysame} {\leavevmode\hbox
to3em{\hrulefill}\thinspace}

\end{document}